
\def\LaTeX{%
  \let\Begin\begin
  \let\End\end
  \let\salta\relax
  \let\finqui\relax
  \let\futuro\relax}

\def\UK{\def\our{our}\let\sz s}
\def\USA{\def\our{or}\let\sz z}

\UK



\LaTeX

\USA


\salta

\documentclass[twoside,a4paper,12pt]{article}
\setlength{\textheight}{24cm}
\setlength{\textwidth}{16cm}
\setlength{\oddsidemargin}{2mm}
\setlength{\evensidemargin}{2mm}
\setlength{\topmargin}{-15mm}
\parskip2mm


\usepackage[usenames,dvipsnames]{color}
\usepackage{amsmath}
\usepackage{amsthm}
\usepackage{amssymb}
\usepackage[mathcal]{euscript}
%
%
\usepackage{cite}
%
%
%


\definecolor{viola}{rgb}{0.3,0,0.7}
\definecolor{ciclamino}{rgb}{1,0,1}
\definecolor{rosso}{rgb}{0.8,0,0}





\def\pcol #1{#1}

\def\revis #1{#1}

\def\gianni #1{#1}
\def\pier #1{#1}
\def\gabri #1{#1}
\def\gult #1{#1}


\bibliographystyle{plain}


%

\finqui

\def\Beq{\Begin{equation}}
\def\Eeq{\End{equation}}
\def\Bsist{\Begin{eqnarray}}
\def\Esist{\End{eqnarray}}

\def\Bthm{\Begin{theorem}}
\def\Ethm{\End{theorem}}
\def\Blem{\Begin{lemma}}
\def\Elem{\End{lemma}}
\def\Bprop{\Begin{proposition}}
\def\Eprop{\End{proposition}}
\def\Bcor{\Begin{corollary}}
\def\Ecor{\End{corollary}}
\def\Brem{\Begin{remark}\rm}
\def\Erem{\End{remark}}

\def\Bdim{\Begin{proof}}
\def\Edim{\End{proof}}
\def\Bcenter{\Begin{center}}
\def\Ecenter{\End{center}}
\let\non\nonumber




\def\step #1 \par{\medskip\noindent{\bf #1.}\quad}


\def\Lip{Lip\-schitz}
\def\Holder{H\"older}
\def\Frechet{Fr\'echet}
\def\aand{\quad\hbox{and}\quad}
\def\lsc{lower semicontinuous}

\def\lhs{left-hand side}
\def\rhs{right-hand side}
\def\sfw{straightforward}


\def\generaliz{generali\sz}
\def\lineariz{lineari\sz}

\def\organiz{organi\sz}


\def\multibold #1{\def\arg{#1}%
  \ifx\arg\pto \let\next\relax
  \else
  \def\next{\expandafter
    \def\csname #1#1#1\endcsname{{\bf #1}}%
    \multibold}%
  \fi \next}

\def\pto{.}

\def\multical #1{\def\arg{#1}%
  \ifx\arg\pto \let\next\relax
  \else
  \def\next{\expandafter
    \def\csname cal#1\endcsname{{\cal #1}}%
    \multical}%
  \fi \next}


\def\multimathop #1 {\def\arg{#1}%
  \ifx\arg\pto \let\next\relax
  \else
  \def\next{\expandafter
    \def\csname #1\endcsname{\mathop{\rm #1}\nolimits}%
    \multimathop}%
  \fi \next}

\multibold
qwertyuiopasdfghjklzxcvbnmQWERTYUIOPASDFGHJKLZXCVBNM.

\multical
QWERTYUIOPASDFGHJKLZXCVBNM.

\multimathop
dist div dom meas sign supp .

\def\calI{\calJ_0}


\def\accorpa #1#2{\eqref{#1}--\eqref{#2}}
\def\Accorpa #1#2 #3 {\gdef #1{\eqref{#2}--\eqref{#3}}%
  \wlog{}\wlog{\string #1 -> #2 - #3}\wlog{}}


\def\separa{\noalign{\allowbreak}}

\def\Neto{\mathrel{{\scriptstyle\nearrow}}}

\def\neto{\mathrel{{\scriptscriptstyle\nearrow}}}
\def\seto{\mathrel{{\scriptscriptstyle\searrow}}}

\def\graffe #1{\mathopen\{#1\mathclose\}}

\def\<#1>{\mathopen\langle #1\mathclose\rangle}
\def\norma #1{\mathopen \| #1\mathclose \|}

\def\iot {\int_0^t}

\def\intQt{\int_{Q_t}}
\def\intQ{\int_Q}
\def\iO{\int_\Omega}
\def\iG{\int_\Gamma}
\def\intS{\int_\Sigma}
\def\intSt{\int_{\Sigma_t}}

\def\dt{\partial_t}
\def\dn{\partial_n}

\def\cpto{\,\cdot\,}

\def\checkmmode #1{\relax\ifmmode\hbox{#1}\else{#1}\fi}
\def\aeO{\checkmmode{a.e.\ in~$\Omega$}}
\def\aeQ{\checkmmode{a.e.\ in~$Q$}}

\def\aeS{\checkmmode{a.e.\ on~$\Sigma$}}
\def\aet{\checkmmode{a.e.\ in~$(0,T)$}}

\def\aaS{\checkmmode{for a.a.~$(t,x)\in\Sigma$}}
\def\aat{\checkmmode{for a.a.~$t\in(0,T)$}}


\def\erre{{\mathbb{R}}}




\def\genspazio #1#2#3#4#5{#1^{#2}(#5,#4;#3)}
\def\spazio #1#2#3{\genspazio {#1}{#2}{#3}T0}
\def\spaziot #1#2#3{\genspazio {#1}{#2}{#3}t0}
\def\spazios #1#2#3{\genspazio {#1}{#2}{#3}s0}
\def\L {\spazio L}
\def\H {\spazio H}
\def\W {\spazio W}
\def\Lt {\spaziot L}
\def\Ht {\spaziot H}

\def\Ls {\spazios L}
\def\C #1#2{C^{#1}([0,T];#2)}


\def\Lx #1{L^{#1}(\Omega)}
\def\Hx #1{H^{#1}(\Omega)}

\def\HxG #1{H^{#1}(\Gamma)}

\def\Cx #1{C^{#1}(\overline\Omega)}

\def\LQ #1{L^{#1}(Q)}
\def\LS #1{L^{#1}(\Sigma)}
\def\CQ #1{C^{#1}(\overline Q)}

\def\Ldue{\Lx 2}
\def\Linfty{\Lx\infty}

\def\Huno{\Hx 1}

\def\Hunoz{{H^1_0(\Omega)}}


\def\LQ #1{L^{#1}(Q)}


\let\theta\vartheta
\let\eps\varepsilon
\let\phi\varphi

\let\TeXchi\chi                         
\newbox\chibox
\setbox0 \hbox{\mathsurround0pt $\TeXchi$}
\setbox\chibox \hbox{\raise\dp0 \box 0 }
\def\chi{\copy\chibox}



\def\suG{{\vrule height 5pt depth 4pt\,}_\Gamma}

\def\thetaG{\theta_\Gamma}
\def\ThetaG{\Theta_\Gamma}
\def\pG{p_\Gamma}
\def\psiG{\psi_\Gamma}

\def\thetaz{\theta_0}
\def\phiz{\phi_0}
\def\thetaQ{\theta_{\!Q}}
\def\phiO{\phi_\Omega}

\def\qz{q_0}

\def\umin{u_{\rm min}}
\def\umax{u_{\rm max}}
\def\thetamin{\theta_\bullet}
\def\thetamax{\theta^\bullet}
\def\phimin{\phi_\bullet}
\def\phimax{\phi^\bullet}

\def\Uad{\calU_{ad}}
\def\kuno{\kappa_1}
\def\kdue{\kappa_2}
\def\uopt{u^*}
\def\thetaopt{\theta^*}
\def\thetaoptG{\theta^*_\Gamma}
\def\phiopt{\phi^*}
\def\xiopt{\xi^*}
\def\ubar{\overline u}
\def\thetabar{\overline\theta}
\def\thetabarG{\thetabar_\Gamma}
\def\phibar{\overline\phi}
\def\Thetabar{\overline\Theta}

\def\thetah{\theta^h}
\def\thetahG{\theta^h_\Gamma}
\def\phih{\phi^h}
\def\zetah{\zeta^h}
\def\zetahG{\zeta^h_\Gamma}
\def\etah{\eta^h}
\def\qbar{\overline q}

\def\ptilde{\tilde p}
\def\qtilde{\tilde q}

\def\vG{v_\Gamma}
\def\wG{w_\Gamma}

\def\un{u_n}
\def\phin{\phi_n}
\def\thetan{\theta_n}
\def\barthetan{\bar\theta_n}
\def\thetanG{\theta_{n,\Gamma}}
\def\xin{\xi_n}

\def\VG{V_\Gamma}
\def\HG{H_\Gamma}
\def\Vp{V^*}

\def\normaV #1{\norma{#1}_V}
\def\normaH #1{\norma{#1}_H}

\let\hat\widehat
\def\Beta{\hat{\vphantom t\smash\beta\mskip2mu}\mskip-1mu}
\def\betaz{\beta^\circ}
\def\betaeps{\beta_\eps}
\def\Betaeps{\Beta_\eps}

\def\thetaeps{\theta_\eps}
\def\thetaepsG{\theta_{\eps,\Gamma}}
\def\phieps{\phi_\eps}
\def\xieps{\xi_\eps}

\def\Pi{\hat\pi}
\let\Lam\Lambda
\def\Lambda{\hat{\vrule height6.5pt depth 0pt width0pt \smash\lambda}}
\def\lambdaeps{\lambda_\eps}
\def\Lambdaeps{\Lambda_\eps}
\def\suplambdaeps{\Lam_\eps}
\def\supLambdaeps{\hat\Lam_\eps}
\def\Reps{R_\eps}

\def\xieps{\xi_\eps}



\Begin{document}


\title{A boundary control problem\\[0.3cm]
  for a possibly singular phase field system\\[0.3cm]
  with dynamic boundary conditions}

\author{}
\date{}
\maketitle
\Bcenter
\vskip-2cm
{\large\sc Pierluigi Colli$^{(1)}$}\\
{\normalsize e-mail: {\tt pierluigi.colli@unipv.it}}\\[.25cm]
{\large\sc Gianni Gilardi$^{(1)}$}\\
{\normalsize e-mail: {\tt \gabri{gianni.gilardi@unipv.it}}}\\[.25cm]
{\large\sc Gabriela Marinoschi$^{(2)}$}\\
{\normalsize e-mail: {\tt gabriela.marinoschi@acad.ro}}\\[.45cm]
$^{(1)}$
{\small Dipartimento di Matematica ``F. Casorati'', Universit\`a di Pavia}\\
{\small \revis{and Research Associate at the IMATI -- C.N.R. PAVIA}}\\
{\small via Ferrata 1, 27100 Pavia, Italy}\\[.2cm]
$^{(2)}$
{\small ``Gheorghe Mihoc-Caius Iacob'' Institute of Mathematical Statistics\\
and \gabri{Applied Mathematics of the Romanian Academy}}\\
{\small Calea 13 Septembrie 13, 050711 Bucharest, Romania}\\[1cm]
\Ecenter

\Begin{abstract}
\pier{This paper deals with \pcol{an optimal control problem related to} a phase field system of Caginalp type with a dynamic boundary condition for the temperature.  \pcol{The control placed in the dynamic boundary condition acts on a part of the boundary.  The analysis carried out in this paper proves the existence of an optimal control for a general class of potentials, possibly singular. The study  includes potentials for which the derivatives may not exist, these being replaced by well-defined subdifferentials. Under some stronger assumptions on the structure parameters and on the potentials (namely for the regular and the logarithmic case having single-valued derivatives), the first order necessary optimality conditions are derived and expressed in terms of the boundary trace of the first adjoint variable.}}

\vskip3mm

\End{abstract}

\salta

\pagestyle{myheadings}
\newcommand\testopari{\sc Colli \ --- \ Gilardi \ --- \ Marinoschi}
\newcommand\testodispari{\sc A boundary control problem for a phase field system}
\markboth{\testodispari}{\testopari}

\finqui


\section{Introduction}
\label{Intro}
\setcounter{equation}{0}

\pier{In this paper we are concerned with an optimal control problem for a nonlinear phase field system of a standard form (cf. the monograph \cite{BrokSpr}), but with a possibly singular double well potential, like in the logarithmic case (cf. the later \eqref{logpot}), and with a dynamic boundary condition for the temperature, in which also the time derivative of the boundary temperature plays a role and where the control variable appears in the external term (see~\eqref{Idynbc}). Let us now introduce and discuss the problem in precise terms.}

A rather general version of the phase field system of \pier{Caginalp type~\cite{Cag}} reads as follows
\Bsist
  & \dt\theta - \Delta\theta + \lambda(\phi) \dt\phi = 0
  & \quad \hbox{in $Q:=(0,T)\times\Omega$}
  \label{caginalpa}
  \\
  & \dt\phi - \sigma \Delta\phi + \calW'(\phi) = \theta \lambda(\phi)
  & \quad \hbox{in $Q$}
  \label{caginalpb}
\Esist
where $\Omega$ is the domain where the evolution takes place,
$T$~is some final time,
$\theta$~denotes the relative temperature around some critical value
that is taken to be $0$ without loss of generality,
and $\phi$ is the order parameter.
Moreover, $\lambda$~is a given real function,
whose meaning is related to the latent heat, and $\sigma$ is a positive constant.
Finally, $\calW'$ represent the derivative of a double-well potential~$\calW$.
Typical example\gabri{s} \pier{for $\calW$} are \pier{the regular potential}
\Beq
 \calW_{reg}(r) = \gianni{r^2(r-1)^2} \,,
  \quad r \in \erre
  \label{regpot}
\Eeq
\pier{with two absolute minima located in $0$ and $1$, and the logarithmic potential}
\Beq   
  \calW_{log}(r) = ((1+r)\ln (1+r)+(1-r)\ln (1-r)) - \pier{a}\, r^2 \,,
  \quad r \in (-1,1)
  \label{logpot}
\Eeq
where the coefficient \pier{$a>0$} is large enough in order to kill convexity.
\gianni{The potential \eqref{regpot} is a \pier{shifted} version of the usual classical potential
given by $r\mapsto\frac 14(r^2-1)^2$ and precisely satisfies our general assumptions given below,
while \eqref{logpot} \pier{has a derivative which behaves singularly in the neighborhoods of $-1$ and $1$}.}
\pier{Generally speaking,} the potential $\calW$ could be just the sum $\calW=\Beta+\Pi$,
where $\Beta$ is a convex function that is allowed to take the value~$+\infty$,
and $\Pi$ is a smooth perturbation (not necessarily concave).
In such a case, $\Beta$~is supposed to be proper and lower semicontinuous
so that its subdifferential \revis{$\beta:=\partial\Beta$} is \revis{well defined} and can replace the derivative
which might not exist.
Of course, \revis{$\calW'$ has to be read as $\beta+\pi$, where $\pi:=\Pi\,{}'$, and}
equation \eqref{caginalpb} becomes a differential inclusion.

Moreover, initial conditions like $\theta(0)=\thetaz$ and $\phi(0)=\phiz$
and suitable boundary conditions must complement the above equations.
As far as the latter are concerned, 
we take the homogeneous Neumann condition for $\phi$ , that~is,
\Beq
  \dn\phi = 0
  \quad \hbox{on $\Sigma := (0,T)\times\Gamma $}
  \label{pier1} 
\Eeq
where $\Gamma$ is the boundary of~$\Omega$
and $\dn$ is the (say, outward) normal derivative. \pier{The position \eqref{pier1} 
is mostly the rule in the literature for the order parameter $\varphi$. Concerning
the temperature $\theta$}, in order to address the boundary control problem,
we choose the following dynamic boundary condition
\Beq
  \dn\theta + \tau \dt\thetaG + \alpha(\thetaG - mu) = 0
  \quad \hbox{on $\Sigma$}
  \label{Idynbc}
\Eeq
where $\thetaG:=\theta\suG$ is the trace of $\theta$ on the boundary
and $u$ is the control, which is supposed to vary in some control box~$\Uad$.
Moreover, in~\eqref{Idynbc}, $\tau$~is a positive time relaxation parameter,
$\alpha$~is a positive constant and $m$ is a nonnegative function defined on~$\Gamma$.
Notice that, in fact, the control $u$ can act just on the subset of $\Gamma$ where $m$ is positive.
Thus, the state system takes the following form
\Bsist
  & \dt\theta - \Delta\theta + \lambda(\phi) \dt\phi = 0
  & \quad \hbox{in $Q$}
  \label{Iprima}
  \\
  & \dt\phi - \sigma \Delta\phi + \beta(\phi) + \pi(\phi) \ni \theta \lambda(\phi)
  & \quad \hbox{in $Q$}
  \label{Iseconda}
  \\
  & \dn\theta + \tau \dt\thetaG + \alpha(\thetaG - mu) = 0
  \aand
  \dn\phi = 0
  & \quad \hbox{on $\Sigma$}
  \label{Ibc}
  \\
  & \theta(0) = \thetaz
  \aand
  \phi(0) = \phiz
  & \quad \hbox{on $\Omega$} .
  \label{Icauchy}
\Esist
\Accorpa\Ipbl Iprima Icauchy
The cost functional we consider depends on two
nonnegative constants $\kuno$~and $\kdue$ and two functions
$\thetaQ$ and $\phiO$ on $Q$ and~$\Omega$, respectively.
We want to minimize
\Beq
  \calJ(u)
  := \frac \kuno 2 \intQ |\theta - \thetaQ|^2
  + \frac \kdue 2 \iO |\phi(T) - \phiO|^2
  \label{Icost}
\Eeq
where $u$, $\theta$ and $\phi$ vary under the constraint of
the state system and $u\in\Uad$,
where the control box $\Uad$ is defined~by
\Beq
  \Uad :=
  \bigl\{ u \in \LS2 : \ \umin\leq u\leq\umax\ \aeS \bigr\}
  \label{Iuad}
\Eeq
for some given bounded functions $\umin$ and $\umax$. \pier{The analysis carried out in this paper shows the existence of an optimal control for a general class of potentials $\calW=\Beta+\Pi$:
indeed, for this purpose $\Beta $ is just assumed to be a general convex and \lsc\ function with minimum $0$ attained at $0$, that is, $\Beta(0)=0$, which is physically reasonable.  On the other hand, the derivation of 
the first order necessary optimality conditions can be made only in case of regular (like~\eqref{regpot}) and singular (like~\eqref{logpot}) potentials. Linearized and adjoint problems are under our investigation and, subsequently, the optimality conditions can be expressed in terms of the adjoint variables (see the Theorem~\ref{CNoptadj} stated in the next section).} 
 
\pier{Let us mention here some related work. As far as we know, the contributions on optimal control problems for phase field models are quite a few and often restricted to the case of regular potentials, or dealing with approximations of the actual systems when the first order optimality conditions are discussed. In this respect, we point out the papers \cite{CoGiMaRo, HoffJiang, HKKY, ZLL} concerned with distributed control problems; we also refer to \cite{BBCG, BCF, CGPS, CGS, CMR, LeSp, SY, SprZheng} for different types of phase field models and other kinds of control problems. The main features of our paper are the study of a boundary control problem and  the consideration of a dynamic boundary condition, in a very simple form: indeed, \eqref{Idynbc} is an affine condition involving the temperature and its time derivative on the boundary, with an external term carrying out the action of the control. About dynamic boundary conditions, also of nonlinear type and possibly involving the Laplace-Beltrami operator, let us quote the articles \cite{BFM,CFS,CF1,CF2,CGS2,CGS3,CS,Ruiz,GrMS}.}

The paper is \organiz ed as follows.
In the next section, we list our assumptions, state the problem in a precise form
and present our results.
The well-posedness of the state system, the regularity results
and the existence of an \pier{optimal} control will be shown
in Sections~\ref{STATE} and~\ref{OPTIMUM}, respectively,
while the rest of the paper is devoted
to the derivation of first order necessary conditions for optimality,
\gianni{which are computed in the case of potentials that \generaliz e
\accorpa{regpot}{logpot}
(some cases of more singular potentials being the subject of a future project of ours)}.
The final result will be proved in Section~\ref{OPTIMALITY}
and it is prepared in Sections~\ref{FRECHET} \pier{with the study of} 
the control-to-state mapping.
Finally, the Appendix is devoted to the rigorous proof
of an estimate that is derived just formally in Section~\ref{STATE}.


\section{Statement of the problem and results}
\label{STATEMENT}
\setcounter{equation}{0}

In this section, we describe the problem under study
and present our results.
As in the Introduction,
$\Omega$~is the body where the evolution takes place.
We assume $\Omega\subset\erre^3$
to~be open, bounded, connected, and smooth,
and we write $|\Omega|$ for its Lebesgue measure.
Moreover, $\Gamma$ and $\dn$ still stand for
the boundary of~$\Omega$ and the outward normal derivative, respectively.
Given a finite final time~$T>0$,
we set for convenience
\Bsist
  && Q_t := (0,t) \times \Omega
  \aand
  \Sigma_t := (0,t) \times \Gamma
  \quad \hbox{for every $t\in(0,T]$}
  \label{defQtSt}
  \\
  && Q := Q_T \,,
  \aand
  \Sigma := \Sigma_T \,.
  \label{defQS}
\Esist
\Accorpa\defQeS defQtSt defQS
Now, we specify the assumptions on the structure of our system.
We assume that
\Bsist
  & \sigma,\, \tau \in (0,+\infty), \quad
  m \in \Lx\infty
  \aand
  m \geq 0
  \quad \aeO
  \label{hplm}
  \\
  & \Beta : \erre \to [0,+\infty]
  \quad \hbox{is convex, proper and l.s.c.}
  \quad \hbox{with} \quad
  \Beta(0) = 0
  \label{hpBeta}
  \\
  & \Pi, \, \Lambda: \erre \to \erre
  \quad \hbox{are \pier{$C^1$ functions}}
  \aand
  {\Pi\,}', \, {\Lambda\,}'
  \quad \hbox{are \Lip\ continuous}
  \label{hpPi}
\Esist
\Accorpa\HPstruttura hplm hpPi
We set for convenience
\Beq
  \beta := \partial\Beta , \quad
  \pi := {\Pi\,}'
  \aand
  \lambda := {\Lambda\,}'
  \label{defbetapi}
\Eeq
and denote by $D(\beta)$ and $D(\Beta)$
the effective domains of $\beta$ and~$\Beta$, respectively.
Moreover, $\betaz(r)$~is the element of $\beta(r)$ having minimum modulus
for every $r\in D(\beta)$
(see, e.g., \cite[p.~28]{Brezis}).
It is well known that $\beta$ is a maximal monotone operator 
\pier{from $\erre$ to $\erre$}
(see, e.g., \cite[Ex.~2.3.4, p.~25]{Brezis}).
Next, in order to simplify notations, we~set\gianni{%
\begin{align}
  V := \Huno, \
  H := \Ldue, \
  W := \graffe{v\in\Hx2: \dn v=0}
  \non 
  \\
   \pier{\hbox{as well as }} \ \HG := L^2(\Gamma)
  \aand
  \VG := H^{1}(\Gamma )
  \label{defspazi}
\end{align}
}%
and endow these spaces with their natural norms.
The symbol $\norma\cpto_X$ stands for the norm in the generic Banach space~$X$,
while $\norma\cpto_p$ is the usual norm in anyone of the $L^p$ spaces
on $\Omega$, $\Gamma$, $Q$ and~$\Sigma$, for $1\leq p\leq\infty$, provided that no confusion can arise.
Furthermore, it is understood that $H$ is embedded in~$\Vp$, the dual space of~$V$,
in the standard way, i.e., in order that
$\<u,v>=\iO uv$ for every $u\in H$ and $v\in V$,
where $\<\cpto,\cpto>$ denotes the duality product between $\Vp$ and~$V$.
Finally, for $v\in\LQ2$ the symbol $1*v$ is the usual \pier{time} convolution,~i.e.,
\Beq
  (1*v)(t) := \iot v(s) \, ds
  \quad \hbox{for $t\in[0,T]$}.
  \label{defconv}
\Eeq

At this point, we describe the state system.
Given $\thetaz$ and $\phiz$ such that
\Beq
  \thetaz \in V , \quad
  \phiz \in W
  \aand
  \betaz(\phiz) \in H
  \label{hpdati}
\Eeq
we look for a quadruplet $(\theta,\thetaG,\phi,\xi)$ satisfying
\Bsist
  & \theta \in \H1H \cap \L\infty V
  \label{regtheta}
  \\
  & \thetaG \in \H1\HG
  \aand
  \thetaG(t) = \theta(t)\suG
  \quad \aat
  \label{regthetaG}
  \\
  & \phi \in \W{1,\infty}H \cap \H1V \cap \L\infty W
  \label{regphi}
  \\
  & \xi \in \L\infty H
  \label{regxi}
\Esist
and solving the problem
\Bsist
  && \iO \dt\theta \, v
  + \iO \nabla\theta \cdot \nabla v
  + \iO \lambda(\phi) \dt\phi \, v
  + \tau \iG \dt\thetaG \, \vG
  + \alpha \iG (\thetaG - mu) \, \vG = 0
  \qquad
  \non
  \\
  && \quad \hbox{for every $(v,\vG)\in V\times\VG$ such that $\vG=v\suG$ and \aet}
  \label{prima}
  \\
  && \dt\phi - \sigma \Delta\phi + \xi + \pi(\phi) = \theta \lambda(\phi)
  \aand
  \xi \in \beta(\phi)
  \quad \aeQ
  \label{seconda}
  \\
  && \dn\phi = 0
  \quad \aeS
  \label{bc}
  \\
  && \theta(0) = \thetaz
  \aand
  \phi(0) = \phiz
  \quad \aeO .
  \label{cauchy}
\Esist
\Accorpa\Regsoluz regtheta regxi
\Accorpa\Pbl prima cauchy
\Accorpa\Tuttopbl regtheta cauchy

\Brem
\label{Strongsoluz}
The variational equation \eqref{prima} is the weak formulation of equation \eqref{Iprima}
and of the dynamic boundary condition contained in~\eqref{Ibc}.
Let us notice that we can deduce both \eqref{Iprima} and the first \pier{condition in}~\eqref{Ibc} from~\eqref{prima}.
Indeed, by writing \eqref{prima} with an arbitrary $v\in\Hunoz$ and $\vG=0$,
we derive \eqref{Iprima} in the sense of distributions on~$Q$.
From~\eqref{regphi} we infer that $\phi$ is bounded since $W\subset\Lx\infty$.
As $\lambda$ is continuous, the same holds for $\lambda(\phi)$,
so that $\lambda(\phi)\dt\phi\in\L2H$.
By comparison in~\eqref{Iprima}, $\Delta\theta$ belongs to~$\L2H$ and \eqref{Iprima} holds \aeQ.
It also follows that the normal derivative $\dn\theta$ makes sense
in a proper Sobolev space of negative order on the boundary
and that the integration--by--parts formula holds in a \generaliz ed sense.
By applying it, one immediately derives the dynamic boundary condition contained in~\eqref{Ibc}
in a \generaliz ed sense.
By comparison in it, $\dn\theta\in\LS2$ and the boundary condition holds \aeS.
Finally, we remark that \eqref{regtheta} and the trace theorem imply $\thetaG\in\L\infty{\HxG{1/2}}$.
\Erem

Our first result ensures well-posedness with the prescribed regularity, stability
and continuous dependence on the control variable in suitable topologies.

\Bthm
\label{Wellposedness}
Assume \HPstruttura\ and \eqref{hpdati}.
Then, for every $u\in\LS2$, problem \Pbl\ has a unique solution $(\theta,\thetaG,\phi,\xi)$
satisfying \Regsoluz,
and the estimate
\Bsist
  && \norma\theta_{\H1H \cap \L\infty V}
  + \norma\thetaG_{\H1\HG}
  \qquad
  \non
  \\
  && \qquad {}
  + \norma\phi_{\W{1,\infty}H \cap \H1V \cap \L\infty W}
  + \norma\xi_{\L\infty H}
  \ \leq \ C_1
  \label{stimasoluz}
\Esist
holds true for some constant $C_1$ that depends only on $\Omega$, $T$,
the structure \HPstruttura\ of the system,
the norms of the initial data associated to~\eqref{hpdati} and an upper bound for~$\norma u_2$.
Moreover, if $u_i\in\LS2$, $i=1,2$, are given
and $(\theta_i,\theta_{i,\Gamma},\phi_i,\xi_i)$ are the corresponding solutions,
then the estimate
\Bsist
  && \norma{\theta_1-\theta_2}_{\L\infty H}
  + \norma{(1*\theta_1)-(1*\theta_2)}_{\L\infty V}
  \non
  \\[1mm]
  && \quad {}
  + \norma{\theta_{1,\Gamma}-\theta_{2,\Gamma}}_{\LS2}
  + \norma{\phi_1-\phi_2}_{\H1H\cap\L\infty V}
  \non
  \\[1mm]
  && \leq C' \, \norma{u_1-u_2}_{\LS2}
  \label{contdep}
\Esist
holds true where $C'$ depends only on $\Omega$, $T$
and on the structure \HPstruttura\ of the system.
\Ethm

\Brem
\label{Phicont}
Since $W$ is compactly embedded in~$\Cx0$,
\revis{the space of continuous functions on~$\overline\Omega$},
the~regularity \eqref{regphi} implies \revis{$\phi\in\CQ0:=\C0{\Cx0}$}
and the estimate $\norma\phi_\infty\leq c\norma\phi_{\L\infty W}$,
where $c$ depends only on~$\Omega$.
Therefore, we also have
\Beq
  \phi \in \CQ0
  \aand
  \norma\phi_\infty \leq C_2
  \label{phicont}
\Eeq
where $C_2$ is a multiple of the constant $C_1$ of~\eqref{stimasoluz}.
\Erem

Some further regularity of the solution is stated in the next result of ours.

\Bthm
\label{Regularity}
Under assumptions \HPstruttura\ and \eqref{hpdati},
the following properties hold true.

$i)$~If in addition $\thetaz\in\Linfty$ and $u\in\LS\infty$, we also have
\Beq
  \theta \in \LQ\infty
  \aand
  \norma\theta_\infty \leq C_3
  \label{thetabdd}
\Eeq
\pier{where $C_3$ is a constant with the same dependencies as $C_1$ and depending} on $\norma\thetaz_\infty$
and on an upper bound for $\norma u_\infty$, in addition.

$ii)$~By also assuming $\betaz(\phiz)\in\Linfty$,
we have that $\xi\in\LQ\infty$ and
\Beq
  \norma\xi_{\LQ\infty} \leq C_4
  \label{stimaxi}
\Eeq
with a constant $C_4$ that depends on the norm
$\norma{\betaz(\phiz)}_\infty$ as well.
\Ethm

The well-posedness result for problem \Pbl\ given by Theorem~\ref{Wellposedness}
allows us to introduce the control-to-state mapping~$\calS$
and to address the corresponding control problem.
We define
\Bsist
  && \calX := \LS\infty
  \aand
  \calY := \calY_1 \times \calY_2 \times \calY_3
  \quad \hbox{where}
  \label{defXY}
  \\
  && \quad \calY_1 := \graffe{v\in\LQ2: \ \ 1*v\in\L\infty V}, \quad
  \calY_2 := \LS2
  \label{defYunodue}
  \\
  && \aand
  \calY_3 := \C0H \cap \L2V
  \label{defYtre}
  \\
  && \calS : \calX \to \calY ,
  \quad
  u \mapsto \calS(u) =: (\theta,\thetaG,\phi)
  \quad \hbox{where}
  \non
  \\
  && \hbox{$(\theta,\thetaG,\phi,\xi)$ is the unique solution to \Tuttopbl\ corresponding to $u$}.
  \qquad
  \label{defS}
\Esist
Next, in order to introduce the control box and the cost functional,
we assume that
\Bsist
  & \umin,\umax\in\LS\infty
  \quad \hbox{satisfy} \quad
  \umin \leq \umax
  \quad \aeS
  \label{hpUad}
  \\
  & \kuno \,,\, \kdue \in [0,+\infty) , \quad
  \thetaQ \in \LQ2
  \aand
  \phiO \in \pier{\Hx1}
  \label{hpJ}
\Esist
\Accorpa\StrutturacCP hpUad hpJ
and define $\Uad$ and $\calJ$
according to the Introduction.
Namely, we~set
\Bsist
  && \Uad :=
  \bigl\{ u \in \calX : \ \umin\leq u\leq\umax\ \aeS \bigr\}
  \phantom \int
  \label{defUad}
  \\
  && \calJ := \calI \circ \calS : \calX \to \erre
  \quad \hbox{where} \quad
  \calI : \calY \to \erre
  \quad \hbox{is defined by}
  \label{defJ}
  \\
  && \calI(\theta,\thetaG,\phi)
  := \frac \kuno 2 \intQ |\theta - \thetaQ|^2
  + \frac \kdue 2 \iO |\phi(T) - \phiO|^2 .
  \label{defI}
\Esist
\Accorpa\Defcontrol defXY defI

Here is our first result on the control problem.

\Bthm
\label{Optimum}
Assume \HPstruttura\ and \eqref{hpdati},
and let $\Uad$ and $\calJ$ be defined \pier{in} \Defcontrol.
Then, there exists $\uopt\in\Uad$ such~that
\Beq
  \calJ(\uopt)
  \leq \calJ(u)
  \quad \hbox{for every $u\in\Uad$} .
  \label{optimum}
\Eeq
\Ethm

From now on, it is understood that the assumptions \HPstruttura\ \pier{and \eqref{hpdati}}
on the structure and on the initial data are satisfied
and that \pier{the} map~$\calS$, the cost functionals $\calI$ and $\calJ$
and the control box $\Uad$ are defined by \Defcontrol.
Thus, we do not remind anything of that in the statements given in the sequel.

Our next aim is to formulate first order necessary optimality conditions.
As $\Uad$ is convex, the desired necessary condition for optimality~is
\Beq
  \< D\calJ(\uopt) , u-\uopt > \geq 0
  \quad \hbox{for every $u\in\Uad$}
  \label{precnopt}
\Eeq
provided that the derivative $D\calJ(\uopt)$ exists in the dual space $\calX^*$
at least in the G\^ateaux sense.
Then, the natural approach consists in proving that
$\calS$ is \Frechet\ differentiable at $\uopt$
and applying the chain rule to $\calJ=\calI\circ\calS$.
We can properly tackle this project under \pier{some} further assumptions
that are satisfied for each of the potentials \accorpa{regpot}{logpot}.
We also have to require something more on~$\lambda$.
Namely, we also suppose that
\begin{gather}
   \hbox{$D(\beta)$ is an open interval and $\beta$ is a single-valued on $D(\beta)$}
  \label{hpDbetabis}
  \\
   \hbox{$\beta$, $\pi$ \pier{and $\lambda$ are $C^2$ functions}}.
  \label{hpnonlinbis}
\end{gather}
\Accorpa\HPstrutturabis hpDbetabis hpnonlinbis
In particular, $\betaz=\beta$.
Furthermore, the \pier{inclusion in} \eqref{seconda} reduces to $\xi=\beta(\phi)$,
and it is no longer necessary to split the nonlinear contribution to the equation
in the form $\xi+\pi(\phi)$.
Hence, we set for brevity
\Beq
  \gamma := \beta + \pi
  \label{defgamma}
\Eeq
\Accorpa\HPstrutturater hpDbetabis defgamma
and observe that $\gamma$ is a $C^2$ function on~$D(\beta)$.

As assumptions \HPstrutturabis\ force $\beta(r)$ to tend to $\pm\infty$
as $r$ tends to a finite end-point of~$D(\beta)$, if any,
we see that combining the further requirement \HPstrutturabis\
with the boundedness of $\phi$ and $\xi$ given by \pier{Theorems~\ref{Wellposedness} and~}\ref{Regularity}
immediately yields

\Bcor
\label{Bddaway}
Under all the assumptions of Theorem~\ref{Regularity},
suppose that \HPstrutturabis\ hold, in addition.
Then, the component $\phi$ of the solution $(\theta,\thetaG,\phi,\xi)$
also satisfies
\Beq
  \phimin \leq \phi \leq \phimax
  \quad \hbox{in $\overline Q$}
  \label{bddaway}
\Eeq
for some constants $\phimin\,,\phimax\in D(\beta)$
that depend only on $\Omega$, $T$,
the structure \HPstruttura\ and \HPstrutturabis\ of the system,
the norms of the initial data associated to \eqref{hpdati},
the norms $\norma\thetaz_\infty$ and $\norma{\beta(\phiz)}_\infty$,
and an upper bound for~$\norma u_\infty$.
\Ecor

As we shall see in Section~\ref{FRECHET},
the computation of the \Frechet\ derivative of $\calS$
leads to the linearized problem that we describe at once
and that can be stated starting from a generic element $\ubar\in\calX$.
Let $\ubar\in\calX$ and $h\in\calX$ be given.
We set $(\thetabar,\thetabarG,\phibar):=\calS(\ubar)$.
Then the linearized problem consists in finding $(\Theta,\ThetaG,\Phi)$
satisfying
\Bsist
  & \Theta \in \H1H \cap \L\infty V 
  \label{regTheta}
  \\
  & \ThetaG \in \H1\HG
  \aand
  \ThetaG(t) = \Theta(t)\suG
  \quad \aat
  \label{regThetaG}
  \\
  & \Phi \in \H1H \cap \L\infty V \cap \L2W
  \label{regPhi}
\Esist
and solving the following problem
\Bsist
  && \iO \dt\Theta \, v
  + \iO \nabla\Theta \cdot \nabla v
  + \iO \lambda(\phibar) \dt\Phi \, v
  + \iO \lambda'(\phibar) \dt\phibar \, \Phi \, v
  \non
  \\
  && \quad {}
  + \tau \iG \dt\ThetaG \, \vG
  + \alpha \iG \ThetaG \, \vG
  = \alpha \iG mh \, \vG
  \non
  \\
  && \quad \hbox{for every $(v,\vG)\in V\times\VG$ such that $\vG=v\suG$ and \aet}
  \label{linprima}
  \\
  && \dt\Phi - \sigma \Delta\Phi + \gamma'(\phibar) \, \Phi
  = \thetabar \lambda'(\phibar) \, \Phi
  + \lambda(\phibar) \, \Theta
  \quad \aeQ
  \label{linseconda}
  \\
  && \dn\Phi = 0
  \quad \aeS
  \label{linbc}
  \\
  && \Theta(0) = 0
  \aand
  \Phi(0) = 0
  \quad \aeO .
  \label{lincauchy}
\Esist
\Accorpa\Reglin regTheta regPhi
\Accorpa\Linpbl linprima lincauchy

\Bprop
\label{Existlin}
Let $\ubar\in\calX$ and $(\thetabar,\thetabarG,\phibar)=\calS(\ubar)$.
Then, for every $h\in\calX$,
there exists a unique triplet $(\Theta,\ThetaG,\Phi)$
satisfying \Reglin\
and solving the linearized problem \Linpbl.
Moreover, the inequality
\Beq
  \norma{(\Theta,\ThetaG,\Phi)}_\calY \leq C_5 \norma h_\calX
  \label{stimaFrechet}
\Eeq
holds true with a constant $C_5$
that depend only on $\Omega$, $T$,
the structure \HPstruttura\ and \HPstrutturabis\ of the system,
the norms of the initial data associated to~\eqref{hpdati},
and the norms
$\norma\ubar_\infty$, $\norma\thetaz_\infty$ and $\norma{\beta(\phiz)}_\infty$.
In particular, the linear map $\calD:h\mapsto(\Theta,\ThetaG,\Phi)$
is continuous from $\calX$ to~$\calY$.
\Eprop

Namely, we shall prove in Section~\ref{FRECHET}
that the map $\calD\in\calL(\calX,\calY)$ introduced in the last \pier{statement} exactly provides the \Frechet\ derivative
$D\calS(\ubar)$ of $\calS$ at~$\ubar$.
Once this is done, we may use the chain rule with $\ubar:=\uopt$
to prove that the necessary condition \eqref{precnopt} for optimality takes the form
\Beq
  \kuno \intQ (\thetaopt - \thetaQ) \Theta
  + \kdue \iO \bigl( \phiopt(T) - \phiO \bigr) \Phi(T)
  \geq 0
  \quad \hbox{for any $u\in\Uad$}
  \qquad
  \label{cnopt}
\Eeq
where $(\thetaopt,\thetaoptG,\phiopt)=\calS(\uopt)$ and, for any given $u\in\Uad$,
the triplet $(\Theta,\ThetaG,\Phi)$ is the solution to the linearized problem
corresponding to $h=u-\uopt$.

The final step then consists in eliminating the pair $(\Theta,\Phi)$ from~\eqref{cnopt}.
This will be done by introducing a triplet $(p,\pG,q)$
that fulfils the regularity requirements
\begin{align}
  & \pier{p \,, q \in \H1H \cap \L\infty V , \quad  q \in  \L2{W}}
  \label{regpq}
  \\
  & \pG \in \H1\HG
  \aand
  \pG(t) = p(t)\suG
  \quad \aat
  \label{regpG}
\end{align}
\Accorpa\Regadj regpq regpG
and solves the following adjoint system:
\Bsist
  && - \iO \dt p \, v
  + \iO \nabla p \cdot \nabla v
  - \iO \lambda(\phiopt) q v
  - \tau \iG \dt\pG \, \vG
  + \alpha \iG \pG \vG
  \non
  \\
  && = \kuno \iO (\thetaopt-\thetaQ) v
  \non
  \\
  && \quad \hbox{for every $v\in V$, $\vG:=v\suG$, and \aet}
  \label{primaadj}
  \\
  && - \iO \dt q \, v
  + \sigma \iO \nabla q \cdot \nabla v
  + \iO \bigl( \gamma'(\phiopt) - \thetaopt \lambda'(\phiopt) \bigr) q v
  = \iO \lambda(\phiopt) \dt p \, v
  \qquad
  \non
  \\
  && \quad \hbox{for every $v\in V$ and \aet}
  \label{secondaadj}
  \\
  && p(T) = 0
  \aand
  q(T) = \kdue \bigl( \phiopt(T) - \phiO \bigr)
  \quad \aeO .
  \label{cauchyadj}
\Esist
\Accorpa\Pbladj primaadj cauchyadj

Clearly, \accorpa{primaadj}{secondaadj} are the variational formulation of a boundary value problem.
Namely, $p$ and $q$ solve two backward parabolic equations
complemented by a dynamic boundary condition for $p$
and the homogeneous Neumann boundary condition for~$q$.
However, it is more convenient to keep the problem in that form.

\Bthm
\label{Existenceadj}
Let $\uopt$ and $(\thetaopt,\thetaoptG,\phiopt)=\calS(\uopt)$
be an optimal control and the corresponding state.
Then the adjoint problem \Pbladj\ has a unique solution $(p,\pG,q)$
satisfying the regularity conditions \Regadj.
\Ethm

We recall that, if $K$ is a closed interval and $y_0\in K$,
the normal cone to $K$ at $y_0$ is the set of $z\in\erre$ such that
$z(y-y_0)\leq0$ for every $y\in K$.
Here is our last result.

\Bthm
\label{CNoptadj}
Let $\uopt$ be an optimal control \pier{and 
$(\thetaopt,\thetaoptG,\phiopt)=\calS(\uopt)$ denote 
the associate state. Moreover, let $(p,\pG,q)$
be} the unique solution to the adjoint problem~\Pbladj\
given by Theorem~\ref{Existenceadj}.
Then, \aaS, we~have
\Beq
  m(x) \pG(t,x) \bigl( y - \uopt(t,x) \bigr) \geq 0
  \quad \hbox{for every $y\in[\umin(t,x),\umax(t,x)]$}
  \label{cnoptadj}
\Eeq
that is, $-m(x)\pG(t,x)$ belongs to the normal cone
to $[\umin(t,x),\umax(t,x)]$ at~$\uopt(t,x)$.
In particular, we have
\Beq
  \gabri{u^{\ast }=u_{\max } , \quad
  u^{\ast }=u_{\min }
  \aand
  u^{\ast }\in [u_{\min },u_{\max }]}
  \non
\Eeq
a.e.\ in the subsets of $\Sigma$
where \revis{$m\pG<0$, $m\pG>0$, $m\pG=0$,} respectively.
\Ethm

In performing our a priori estimates in the remainder of the paper,
we often use the \Holder\ inequality
(with the standard notation $p'$ for the conjugate exponent of~$p$),
its consequences and the elementary Young inequalities
\begin{align}
  ab \leq \gianni\omega \, a^{1/\gianni\omega} + (1-\gianni\omega) \, b^{1/(1-\gianni\omega)}
  \aand
  ab \leq \delta a^2 + \frac 1 {4\delta} \, b^2
  \non
  \\
  \quad \hbox{for every $a,b\geq 0$, \ $\gianni\omega\in(0,1)$ \ and \ $\delta>0$}
  \label{young}
\end{align}
as well as the continuous (in fact compact) embedding $V\subset\Lx4$.
Moreover, in order to avoid a boring notation,
we follow a general rule to denote constants.
The small-case symbol $c$ stands for different constants which depend only
on~$\Omega$, on the final time~$T$, the shape of the nonlinearities
and on the constants and the norms of
the functions involved in the assumptions of our statements.
A~small-case symbol with a subscript like $c_\delta$
indicates that the constant might depend on the parameter~$\delta$, in addition.
Hence, the meaning of $c$ and $c_\delta$ might
change from line to line and even in the same chain of equalities or inequalities.
On the contrary, we mark precise constants which we can refer~to
by using different symbols, e.g., capital letters.


\section{The state system}
\label{STATE}
\setcounter{equation}{0}

This section is devoted to the \pier{proofs} of Theorems~\ref{Wellposedness} and~\ref{Regularity}.
As far as the former is concerned,
we start proving its second part, i.e., the continuous dependence formula \eqref{contdep}.
From this we derive uniqueness as well.

\step
Continuous dependence and uniqueness

We first derive an identity that is satisfied by any solution.
By integrating \eqref{prima} with respect to time, we obtain
\Bsist
  && \iO \theta v
  + \iO \nabla(1*\theta) \cdot \nabla v
  + \tau \iG \thetaG \vG
  + \alpha \iG \bigl( 1*(\thetaG-mu) \bigr) \vG
  \non
  \\
  \separa
  && = \iO \thetaz \, v
  + \tau \iG (\thetaz\suG) \, \vG
  - \iO \bigl( \Lambda(\phi) - \Lambda(\phiz) \bigr) v
  \non
  \\
  && \quad \hbox{for every $(v,\vG)\in V\times\VG$ such that $\vG=v\suG$ and \aet}.
  \label{intprima}
\Esist
Now, we pick two elements $u_i\in\LS2$, $i=1,2$,
and consider two corresponding solutions
$(\theta_i,\theta_{i,\Gamma},\phi_i,\xi_i)$.
We write \eqref{intprima} for both controls and solutions
and test the difference by choosing
$v=\theta:=\theta_1-\theta_2$ and $\vG=\thetaG:=\theta_{1,\Gamma}-\theta_{2,\Gamma}$.
Then, we integrate over $(0,t)$, where $t\in(0,T)$ is arbitrary.
At the same time, we write \eqref{seconda} for both solutions,
multiply the difference by $\phi:=\phi_1-\phi_2$
and integrate over~$Q_t$.
Finally, we add the equalit\pier{ies} we have obtained to each other.
By also setting $u:=u_1-u_2$ and $\xi:=\xi_1-\xi_2$ for brevity, we \pier{infer~that}
\Bsist
  && \intQt |\theta|^2
  + \frac 12 \iO |\pier{\nabla (1*\theta)}(t)|^2
  + \frac \tau 2 \intSt |\thetaG|^2
  + \frac \alpha 2 \iG |(1*\thetaG)(t)|^2
  \non
  \\
  && \quad {}
  + \frac 12 \iO |\phi(t)|^2
  + \sigma \intQt |\nabla\phi|^2
  + \intQt \xi \phi
  \non
  \\
  \separa
  && = \alpha \intSt m (1*u) \thetaG
  - \intQt \bigl( \Lambda(\phi_1) - \Lambda(\phi_2) \bigr) \theta
  \non
  \\
  && \quad {}
  - \intQt \bigl( \pi(\phi_1) - \pi(\phi_2) \bigr) \phi
  + \intQt \bigl( \theta_1 \lambda(\phi_1) - \theta_2 \lambda(\phi_2) \bigr) \phi .
  \label{percontdep}
\Esist
All the terms on the \lhs\ are nonnegative, \pier{including} the last one since $\beta$ is monotone.
We estimate each term on the \rhs, separately.
In the sequel, $\delta$~is a positive parameter.
We~have
\Bsist
  && \alpha \intSt m (1*u) \thetaG
  = - \alpha \intSt mu (1*\thetaG)
  + \alpha \iG m (1*u)(t) \, (1*\thetaG)(t)
  \non
  \\
  && \leq c \intSt |1*\thetaG|^2
  + c \intSt |u|^2
  + \delta \iG |(1*\thetaG)(t)|^2
  + c_\delta \iG |(1*u)(t)|^2
  \label{perraggiopalla}
\Esist
and the last integral is bounded by $\pier{c_\delta} \norma u^2_{\LS2}$
due to the \Holder\ inequality.
Next, \pier{owing} to the boundedness of $\phi_1$ and $\phi_2$
ensured by Remark~\ref{Phicont} and to the regularity of $\Lambda$ on bounded intervals,
we infer~that
\Bsist
  && - \intQt \bigl( \Lambda(\phi_1) - \Lambda(\phi_2) \bigr) \theta
  \leq c \intQt |\phi| \, |\theta|
  \leq \delta \intQt |\theta|^2
  + c_\delta \intQt |\phi|^2
  \non
\Esist
and the third integral on the \rhs\ of \eqref{percontdep}
can be treated in a similar way.
Finally, we deal with the last term.
As $\lambda$ is \Lip\ continuous (see~\eqref{hpPi})
and $\phi_2$ is bounded,
we~have
\Bsist
  && \intQt \bigl( \theta_1 \lambda(\phi_1) - \theta_2 \lambda(\phi_2) \bigr) \phi
  = \intQt \theta_1 \bigl( \lambda(\phi_1) - \lambda(\phi_2) \bigr) \phi
  + \intQt \theta \lambda(\phi_2) \phi
  \non
  \\
  && \leq c \intQt |\theta_1| \, |\phi|^2
  + c \intQt |\theta| \, |\phi|
  \leq c \intQt |\theta_1| \, |\phi|^2
  + \delta \intQt |\theta|^2
  + c_\delta \intQt |\phi|^2 .
  \non
\Esist
On the other hand, by combining the \Holder\ inequality,
the Sobolev inequality $\norma v_4\leq c\normaV v$ for every $v\in V$,
and the regularity $\theta_1\in\L\infty V$,
we obtain
\Bsist
  && \intQt |\theta_1| \, |\phi|^2
  \leq \iot \norma{\theta_1(s)}_4 \, \norma{\phi(s)}_4 \, \norma{\phi(s)}_2 \, ds
  \non
  \\
  && \leq c \iot \norma{\phi(s)}_V \, \norma{\phi(s)}_2 \, ds
  \leq \delta \intQt |\nabla\phi|^2
  + c_\delta \intQt |\phi|^2 .
  \non
\Esist
At this point, we collect all the estimates we have derived,
choose $\delta$ small enough and apply the Gronwall lemma.
\pier{Thus, we} obtain~\eqref{contdep} and uniqueness easily follows.
Indeed, taking $u_1=u_2$ in \eqref{contdep} immediately yields
$\theta_1=\theta_2$, $\theta_{1,\Gamma}=\theta_{2,\Gamma}$ and $\phi_1=\phi_2$.
By comparison in~\eqref{seconda},
we also have $\xi_1=\xi_2$.\qed

\medskip

In order to complete the proof of Theorem~\ref{Wellposedness},
we have to show the existence of a solution
and to establish estimate~\eqref{stimasoluz}.
To this end, we first replace $\beta$, $\Lambda$ and $\lambda$ by the smooth approximation of them
$\betaeps$, $\Lambdaeps$ and $\lambdaeps$ we introduce below,
where $\eps$ is a positive parameter, say $\eps\in(0,1)$.
By doing that, we obtain the approximating problem of finding
a quadruplet $(\thetaeps,\thetaepsG,\phieps,\xieps)$ satisfying
regularity requirements of type \Regsoluz\
and solving
\Bsist
  && \iO \dt\thetaeps \, v
  + \iO \nabla\thetaeps \cdot \nabla v
  + \iO \dt\Lambdaeps(\phieps) \, v
  + \tau \iG \dt\thetaepsG \, \vG
  + \alpha \iG (\thetaepsG - mu) \, \vG = 0
  \qquad
  \non
  \\
  && \quad \hbox{for every $(v,\vG)\in V\times\VG$ such that $\vG=v\suG$ and \aet}
  \label{primaeps}
  \\
  && \dt\phieps - \sigma \Delta\phieps + \xieps + \pi(\phieps) = \thetaeps \lambdaeps(\phieps)
  \aand
  \xieps = \betaeps(\phieps)
  \quad \aeQ
  \label{secondaeps}
  \\
  && \dn\phieps = 0
  \quad \aeS
  \label{bceps}
  \\
  && \thetaeps(0) = \thetaz
  \aand
  \phieps(0) = \phiz
  \quad \aeO  .
  \label{cauchyeps}
\Esist
\Accorpa\Pbleps primaeps cauchyeps
The regularity we require for the solution is the following
\Bsist
  & \thetaeps \in \H1H \cap \L\infty V
  \label{regthetaeps}
  \\
  & \thetaepsG \in \H1\HG
  \aand
  \thetaepsG(t) = \thetaeps(t)\suG
  \quad \aat
  \label{regthetaepsG}
  \\
  & \phieps \in \H1H \cap \L\infty V \cap \L2W .
  \label{regphieps}
\Esist
The lower level \eqref{regphieps} with respect to \eqref{regphi}
has been chosen for convenience.
However, the solution we find also satisfies \eqref{regphi},
as it will be clear from the proof.
In the above equations,
$\betaeps$~is the Yosida regularization of $\beta$ at level~$\eps$
(see, e.g., \cite[p.~28]{Brezis}).
It is well known that
$\betaeps$ is maximal monotone, single-valued and \Lip\ \pier{continuous}.
We also introduce the function $\Betaeps$ defined~by
\Beq
  \Betaeps(r) := \int_0^r \betaeps(s) \, ds
  \quad \hbox{for $r\in\erre$}
  \label{defBetaeps}
\Eeq
and recall that
\Bsist
  && |\betaeps(r)| \leq |\betaz(r)|
  \aand
  \betaeps(r)\to\betaz(r)
  \quad \hbox{for $r\in D(\beta)$}
  \label{propYosida}
  \\
  \noalign{\smallskip}
  && 0 \leq \Betaeps(r) \leq \Beta(r)
  \quad  \hbox{for every $r\in\erre$}.
  \label{propBetaeps}
\Esist
For \eqref{propYosida} see, e.g., \cite[Prop.~2.6, p.~28]{Brezis},
while \eqref{propBetaeps} follows from \eqref{propYosida} and $\betaeps(0)=0$ (\pier{cf.}~\eqref{hpBeta}).
Furthermore, $\Lambdaeps$ as well as its derivative $\lambdaeps$ are defined~by
\begin{align}
  & \Lambdaeps(r) := \Lambda(r) \, \zeta(\eps r)
  \aand
  \lambdaeps(r) := \frac d{dr} \, \Lambdaeps(r)
  \quad \hbox{for $r\in\erre$,\quad where}
  \non
  \\
  & \quad \hbox{$\zeta\in C^{\gult{\infty}}(\erre)$ \ satisfies\quad
      $\zeta(r)=1$ for $|r|<1$ and $\zeta(r)=0$ for $|r|>2$}.
  \label{deflambdaeps}
\end{align}
Notice that both $\Lambdaeps$ and $\lambdaeps$ are bounded and \Lip\ continuous,
and we~set
\Beq
  \supLambdaeps := \sup |\Lambdaeps|
  \aand
  \suplambdaeps := \sup |\lambdaeps| .
  \label{suplambda}
\Eeq

Our project is the following:
$i)$~we prove that problem \Pbleps\ has at least a solution
by a \pier{fixed} point argument;
$ii)$~\pier{using} compactness and monotonicity methods
we \pier{show} that its solution tends to a solution of problem \Pbl\ as $\eps\seto0$,
at least for a subsequence.
We need two lemmas.

\Blem
\label{Thetaphi}
Let $\thetaeps\in\L2H$.
Then, there exists a unique $\phieps$ satisfying \eqref{regphieps}, \pier{\eqref{secondaeps}--\eqref{bceps}}
and the second Cauchy condition \pier{in~\eqref{cauchyeps}.}
Moreover, the estimate
\Beq
  \norma\phieps_{\H1H\cap\L\infty V\cap\L2W}
  \leq C_\eps \bigl( 1 + \norma\thetaeps_{\L2H} \bigr)
  \label{thetaphi}
\Eeq
holds true with a constant $C_\eps$ that depends
on the structure \HPstruttura, the norms involved in \eqref{hpdati}
and $\eps$, but it is independent of~$\thetaeps$.
\Elem

\Bdim
We are dealing with a standard \pier{semilinear} parabolic problem
that has a unique solution with the required regularity.
We just derive estimate \eqref{thetaphi} and control the dependence of constants.
We multiply equation \eqref{secondaeps} by $\dt\phieps$ 
\gianni{and add the same integral to both sides, for convenience.
We~have}
\begin{align}
  & \intQt |\dt\phieps|^2 \pier{{}+\intQt |\phieps|^2}
  + \frac \sigma 2 \iO |\nabla\phieps(t)|^2
  + \iO \Betaeps(\phieps(t))
  \non
  \\
  & = \frac \sigma 2 \iO |\nabla\phiz|^2
  + \iO \Betaeps(\phiz)
  - \intQt \pi(\phieps) \dt\phieps
  + \intQt \thetaeps \, \lambdaeps(\phieps) \, \dt\phieps
  \gianni{{}+\intQt |\phieps|^2}
  \label{perthetaphi}
\end{align}
where $\Betaeps$ is given by~\eqref{defBetaeps}.
\revis{We recall~\eqref{propBetaeps} and \eqref{hpBeta}, which imply that 
$\Betaeps(\phiz)\leq\Beta(\phiz)\leq \beta^\circ(\phiz)\phiz$ a.e.\ in~$\Omega$. 
Thus \eqref{hpdati} yields $\norma{\Betaeps(\phiz)}_1\leq c$.
Furthermore, notice}
that $|\lambdaeps(\phieps)|\leq\suplambdaeps$ (see \eqref{suplambda}).
On the other hand, $\pi$~has a linear growth, so that
\Beq
  - \intQt \pi(\phieps) \dt\phieps
  \leq \frac 12 \intQt |\dt\phieps|^2 + c \intQt |\phieps|^2 + c 
  \non
\Eeq
and the last integral, which coincides with the last term on the \rhs\ of \eqref{perthetaphi},
can be treated \pier{as follows: for} every $s\in[0,t]$, we~have
\Bsist
  && \phieps(s) = \phiz + \int_0^s \dt\phieps(s') \, ds' , 
  \quad \hbox{whence}
  \non
  \\
  && |\phieps(s)|^2
  \leq 2 |\phiz|^2 + 2 \, \Bigl| \int_0^s \dt\phieps(s') \, ds' \Bigr|^2
  \leq c + c \int_0^s |\dt\phieps(s')|^2 \, ds' \,.
  \non
\Esist
It \pier{turns out} that
\Beq
  \iO |\phieps(t)|^2 \leq c + c \intQt |\dt\phieps|^2
  \aand
  \intQt |\phieps|^2 
  \leq c + c \iot \Bigl( \int_{Q_s} |\revis\dt\phieps|^2 \Bigr) \, ds \,.
  \label{auxil}
\Eeq
The second \pier{inequality in}~\eqref{auxil} implies that its \lhs,
i.e., the integral we are dealing with, can be handled by the Gronwall lemma
and we conclude that
\Beq
  \norma\phieps_{\H1H}
  + \norma{\nabla\phieps}_{\L\infty H}
  \leq c + c \bigl( 1 + \suplambdaeps \bigr) \norma\thetaeps_{\L2H}
  \non
\Eeq
where $c$ depends only on the structure and the initial datum~\pier{$\phiz$}.
\gianni{From this estimate and the first \pier{inequality in}~\eqref{auxil} it} follows that
\Beq
  \norma\phieps_{\H1H\cap\L\infty V}
  \leq c + c \bigl( 1 + \suplambdaeps \bigr) \norma\thetaeps_{\L2H}
  \non
\Eeq
with a similar new constant~$c$.
Now, \aat, we write \eqref{secondaeps} at time~$t$,
multiply by $\xieps(t)$ and integrate over~$\Omega$.
We obtain \aet
\Bsist
  && \revis\sigma \iO \betaeps'(\phieps) |\nabla\phieps|^2
  + \iO |\xieps|^2
  = \iO \bigl( - \dt\phieps - \pi(\phieps) + \thetaeps \, \lambdaeps(\phieps) \bigr) \, \xieps
  \non
  \\
  && \leq \frac 12 \iO |\xieps|^2
  + c \, \iO \bigl( |\dt\phieps|^2 + 1 + |\phieps|^2 + |\thetaeps\lambdaeps(\phieps)|^2 \bigr)
  \non
\Esist
whence immediately \aat
\Beq
  \normaH{\xieps(t)}
  \leq c \bigl(
    1
    + \normaH{\dt\phieps(t)}
    + \normaH{\phieps(t)}
    + \normaH{\thetaeps(t)\lambdaeps(\phieps(t))}
  \bigr) .
  \label{stimaxipuntuale}
\Eeq
\pier{In view of} the regularity of $\phieps$ already achieved and the one of~$\thetaeps${\pier ,}
recalling \eqref{suplambda} we deduce~that
\Beq
  \norma\xieps_{\L2H}
  \leq c + c \, (1 + \suplambdaeps) \, \norma\thetaeps_{\L2H}
  \non
\Eeq
where $c$ has the same dependenc\pier{ies} as above.
As the estimate of the norm of $\phieps$ in $\L2W$ follows
by comparison in~\eqref{secondaeps} and elliptic regularity,
the proof of \eqref{thetaphi} is complete.
\Edim

\Blem
\label{Phitheta}
Let $\phieps\in\H1H$.
Then, there exists a unique pair $(\thetaeps,\thetaepsG)$ satisfying
\accorpa{regthetaeps}{regthetaepsG} and the first initial condition \pier{in}~\eqref{cauchyeps}
and solving the variational equation~\eqref{primaeps}.
Moreover, the estimates
\begin{gather}
 \norma\thetaeps_{\LQ2} \leq \Reps
  \label{raggiopalla}
  \\
  \norma\thetaeps_{\H1H\cap\L\infty V}
  + \norma\thetaepsG_{\H1\HG}
  \leq D_\eps \bigl( 1 + \norma\phieps_{\H1H} \bigr)
  \label{phitheta}
\end{gather}
hold true with constant $\Reps$ and $D_\eps$ that depend
on the structure \HPstruttura, the norms involved in \eqref{hpdati},
the norm $\norma u_2$ and $\eps$, but they are independent of~$\phieps$.
\Elem

\Bdim
We set
\Beq
  \calV := \graffe{(v,\vG) \in V \times \revis\HG :\ \vG=v\suG} ,
  \quad
  \calH := H \times \HG
  \label{calVH}
\Eeq
and endow \pier{these} spaces with the scalar products defined on $\calV^2$ and $\calH^2$ by
\Bsist
  && \bigl( (w,\wG) , (v,\vG) \bigr)_{\calV}
  := \iO (\nabla w \cdot \nabla v + wv) + \iG \wG \vG
  \label{prodV}
  \\
  && \bigl( (w,\wG) , (v,\vG) \bigr)_{\calH}
  := \iO wv + \tau \iG \wG \vG
  \label{prodH}
\Esist
respectively.
Then, we obtain two Hilbert spaces and $\calV$ is continuously and densely embedded in~$\calH$,
so that we can construct the Hilbert triplet $(\calV,\calH,\calV^*)$ in the usual way.
Moreover, we define the continuous bilinear form $a$ on $\calV^2$ by the formula
\Beq
  a \bigl( (w,\wG) , (v,\vG) \bigr)
  := \iO \nabla w \cdot \nabla v + \alpha \iG \wG \vG
  \label{defforma}
\Eeq
and consider the operator $\calA\in\calL(\calV,\calV^*)$ associated to~$a$.
As $a$ clearly satisfies
\Beq
  a \bigl( (v,\vG) , (v,\vG) \bigr)
  + \bigl( 1 + {\textstyle \frac 1\tau} \bigr) \norma{(v,\vG)}_{\calH}^2
  \geq \norma{(v,\vG)}_{\calV}^2
  \quad \hbox{for every $(v,\vG)\in\calV$}
  \non
\Eeq
the general theory \pier{(see, e.g., \cite{Lions})} ensures that,
for every $F\in\L2{\calV^*}$ and $U_0\in\calH$,
there exists a unique $U$ satisfying
\Bsist
  && U \in \H1{\calV^*} \cap \L2\calV \subset \C0\calH
  \non
  \\
  && U'(t) + \calA \, U(t) = F(t)
  \quad \aet
  \aand
  U(0) = U_0
  \non
\Esist
and that $U\in\H1\calH\cap\L\infty\calV$
whenever $F\in\L2\calH$ and $U_0\in\calV$.
By accounting for our assumption on $\phieps$ and~\eqref{hpdati},
we choose
\Beq
  F = \bigl( \gianni -\dt\Lambdaeps(\phieps) , \alpha m u \bigr)
  \in \L2\calH
  \aand
  U_0 = \bigl( \thetaz \,, \thetaz\suG \bigr)
  \in \calV
  \non
\Eeq
and obtain a unique pair $(\thetaeps,\thetaepsG)$ satisfying
\accorpa{regthetaeps}{regthetaepsG} and the first initial condition~\eqref{cauchyeps}
and solving the variational equation~\eqref{primaeps}.
Let us now prove estimates \accorpa{raggiopalla}{phitheta}.
We observe that the \pier{analog} of \eqref{intprima}
obtained by replacing $\Lambda$ by $\Lambdaeps$
holds for $(\thetaeps,\thetaepsG,\phieps)$.
So, we \pier{take~$(\thetaeps,\thetaepsG)$ as test pair} and have
\begin{align}
  & \intQt |\thetaeps|^2
  + \frac 12 \iO |\nabla(1*\thetaeps)(t)|^2
  + \tau \intSt |\thetaepsG|^2
  + \frac \alpha 2 \iG |(1*\thetaepsG)(t)|^2
  \non
  \\
  & = \alpha \intSt m (1*u) \thetaepsG
  + {}\gult{\intQt \thetaz \thetaeps
  + \tau \intSt (\thetaz\suG) \thetaepsG
  + \intQt (\Lambdaeps(\phiz) - 
  \Lambdaeps(\phieps) ) \thetaeps} \,.
  \label{pier2}
\end{align}
All the integrals on the \lhs\ are nonnegative
and the first term on the \rhs\ can be estimated as we did in~\eqref{perraggiopalla}, namely
\Beq
  \alpha \intSt m (1*u) \thetaepsG
  \leq c \intSt |1*\thetaepsG|^2
  + c \intSt |u|^2
  + \delta \iG |(1*\thetaepsG)(t)|^2
  + c_\delta \iG |(1*u)(t)|^2
  \non
\Eeq
and the last integral is bounded by $\pier{c_\delta} \norma u^2_{\LS2}$
due to the \Holder\ inequality.
The remaining terms on the \rhs\ \gult{of \eqref{pier2}} are treated in the usual way
\pier{with the help of} the \Holder\ and Young inequalities\gult{; just 
for the last term we point out that 
$$
\intQt (\Lambdaeps(\phiz) - \Lambdaeps(\phieps) ) \thetaeps
\leq 2 \supLambdaeps \intQt |\thetaeps| \leq \delta \intQt |\thetaeps|^2
+ c_\delta \big|\supLambdaeps\big|^2 
$$
thanks to \eqref{suplambda}.}
Thus, by choosing $\delta$ small enough and applying the Gronwall lemma,
\gult{it is straightforward to} obtain the desired estimate~\eqref{raggiopalla}.
Let us now prove~\eqref{phitheta}.
The rigorous argument could rely on testing \eqref{primaeps}
by a $\calV$-valued approximation of $(\dt\thetaeps,\dt\thetaepsG)$.
However, we prefer to avoid such a detail and formally test the equation
by $(\dt\thetaeps,\dt\thetaepsG)$, directly.
We~have
\Bsist
  && \intQt |\dt\thetaeps|^2
  + \frac 12 \iO |\nabla\thetaeps(t)|^2
  + \tau \intSt |\dt\thetaepsG|^2
  + \frac \alpha 2 \iG |\thetaepsG(t)|^2
  \non
  \\
  && = - \intQt \lambdaeps(\phieps) \dt\phieps \, \dt\thetaeps
  + \alpha \intSt m u \, \dt\thetaepsG
  + \frac 12 \iO |\nabla\thetaz|^2
  + \frac \alpha 2 \iG |\thetaz\suG|^2 .
  \non
\Esist
As $|\lambdaeps(\phieps)|\leq\suplambdaeps$ (see \eqref{suplambda})
and $\thetaz$ satisfies~\eqref{hpdati}, we immediately derive~\eqref{phitheta}
\pier{owing} to the \Holder\ and Young inequalities.
\Edim

At this point, we are ready to complete the proof of Theorem~\ref{Wellposedness}
by following our project sketched above.

\step
Existence of the approximating solution

As said before, we are going to use a \pier{fixed} point argument.
We often avoid stressing the d\gabri{e}pendence on~$\eps$, which is fixed.
We consider the closed ball of $\pier{\LQ2}$
\Beq
  \calB := \graffe{v\in\pier{\LQ2}:\ \norma v_2\leq\Reps}
  \label{defpalla}
\Eeq
where $\Reps$ is given by Lemma~\ref{Phitheta} (see \eqref{raggiopalla})
and define the map $\calF:\calB\to\pier{\LQ2}$ by the following steps:
$i)$~for $\bar\theta\in\pier{\calB}$ we apply Lemma~\ref{Thetaphi}
where $\thetaeps$ is replaced by~$\bar\theta$,
find the solution~$\phieps$ and term it $\pier{\calE}(\bar\theta)$;
$ii)$~by starting from such a~$\phieps$, we apply Lemma~\ref{Phitheta},
find the solution $(\thetaeps,\thetaepsG)$
and set $\calF(\bar\theta):=\thetaeps$.
By construction, \gult{it turns out that}  $\calF(\bar\theta)\in\calB$\gult{: indeed, the constant $\Reps$ in \eqref{raggiopalla} is independent of $\norma{\phieps}_{\H1H}$.} 
Moreover, with the above notation, we deduce from Lemmas~\ref{Thetaphi} and~\ref{Phitheta} \gult{(cf.~\eqref{thetaphi} and \eqref{phitheta})}
\begin{gather}
  \norma\phieps_{\H1H\cap\L\infty V\cap\L2W}
  \leq C_\eps \bigl( 1 + \Reps \bigr)
  \label{stimaphieps}
  \\
  \norma\thetaeps_{\H1H\cap\L\infty V}
  + \norma\thetaepsG_{\H1\HG}
  \leq D_\eps \bigl( 1 + C_\eps (1 + \Reps) \bigr).
  \label{stimathetaeps}
\end{gather}
Therefore, $\calF(\calB)$ is relatively compact by the Aubin--Lions lemma
(see, e.g., \cite[Thm.~5.1, p.~58]{Lions2}).
Now, we verify that $\calF$ is continuous.
So, we assume that $\barthetan\in\calB$
and that $\barthetan$ converges to $\bar\theta$ in $\LQ2$,
and we prove that $\thetan:=\calF(\barthetan)$ converges to $\calF(\bar\theta)$ in~\revis{$\LQ2$}.
We set for convenience $\phin:=\pier{\calE}(\thetan)$ and $\thetanG:=\thetan\suG$.
As $\barthetan\in\calB$, estimates \accorpa{stimaphieps}{stimathetaeps}
hold for $\phin$, $\thetan$ and~$\thetanG$.
Hence, for a subsequence, $\phin$, $\thetan$ and $\thetanG$
converge to some $\phi$, $\theta$ and $\thetaG$
in the correponding weak or weak star topologies.
Clearly, $\thetaG=\theta\suG$.
Moreover, $\phin$ and $\thetan$ converge to $\phi$ and $\theta$ strongly in $\LQ2$
by the Aubin--Lions lemma.
This implies that $f(\phin)$ converges to $f(\phi)$ strongly in $\LQ2$
for any \Lip\ continuous function~$f:\erre\to\erre$,
and this is the case if $f$ is either
$\lambdaeps$, or $\Lambdaeps$, or $\betaeps$, or~$\pi$.
As a consequence, $\dt\Lambdaeps(\phin)=\lambdaeps(\phin)\dt\phin$ and
$\thetan\lambdaeps(\phin)$ converge to $\pier{\lambdaeps}(\phi)\dt\phi$ and $\theta\lambdaeps(\phi)$
at least weakly in~$\LQ1$.
Therefore, the quadruplet $(\theta,\thetaG,\phi,\xi)$, where $\xi:=\betaeps(\phi)$,
solves the integrated version of the variational formulation
of the approximating problem \Pbleps\ with smooth test functions,
with $\thetaeps$ replaced by $\bar\theta$ in \eqref{secondaeps}.
This easily implies that $\theta=\calF(\bar\theta)$.
As the same argument holds for any subsequen\gabri{c}e
extracted from~$\graffe{\thetan}$, the continuity we have claimed is proved.
Therefore, we can apply the Schauder \pier{fixed} point theorem
and conclude that $\calF$ has a \pier{fixed} point~$\thetaeps$.
As $\thetaeps$ belongs to $\calF(\calB)$,
it satisfies~\eqref{regthetaeps}.
Moreover, $\thetaepsG:=\thetaeps\suG$ belongs to $\H1\HG$,
the function $\phieps:=\pier{\calE}(\thetaeps)$ satisfies~\eqref{regphieps}
and the quadruplet $(\thetaeps,\thetaepsG,\phieps,\xieps)$,
where $\xieps:=\betaeps(\phieps)$, solves \Pbleps.\qed

The last step consists in letting $\eps$ tend to zero and getting a solution to \Pbl.
To this aim, we derive a~priori estimates on the approximating solution
that are uniform with respect to~$\eps\in(0,1)$.
As such estimates are conserved in the limit as $\eps\seto0$,
inequality \eqref{contdep} is established as a consequence.

\step
First a priori estimate

We choose $v=\thetaeps$ and $\vG=\thetaepsG$ in \eqref{primaeps} and integrate over~$(0,t)$.
At the same time, we multiply \eqref{secondaeps} by $\dt\phieps$ and integrate over~$Q_t$.
Then, we add the equalities obtained \pier{in} this way and observe that
the terms involving the nonlinear function $\lambdaeps$ cancel out.
Hence, we~have
\Bsist
  && \frac 12 \iO |\thetaeps(t)|^2
  + \intQt |\nabla\thetaeps|^2
  + \frac \tau 2 \iG |\thetaepsG(t)|^2
  + \alpha \intSt |\thetaepsG|^2
  \non
  \\
  && \quad {}
  + \intQt |\dt\phieps|^2
  + \frac \sigma 2 \iO |\nabla\phieps(t)|^2
  + \iO \Betaeps(\phieps(t))
  \non
  \\
  \separa
  && = \frac 12 \iO |\thetaz|^2
  + \frac \tau 2 \iG |\thetaz\suG|^2
  + \frac \sigma 2 \iO |\nabla\phiz|^2
  + \iO \Betaeps(\phiz) 
  \non
  \\
  && \quad {}
  \gianni{\pier{{}+ \alpha \intSt mu \thetaepsG} {} -\intQt \pi(\phieps)\dt\phieps
  \,.}
  \non
\Esist
As \eqref{propBetaeps} holds,
\gianni{$\pi$ is \Lip\ continuous, $m$ and $u$ are bounded},
and the data satisfy~\eqref{hpdati},
we \gianni{easily} deduce~that
\Beq
  \norma\thetaeps_{\L\infty H\cap\L2V}
  + \norma\thetaepsG_{\L\infty\HG}
  + \norma\phieps_{\H1H\cap\L\infty V}
  \leq c \,.
  \label{primastima}
\Eeq
Since $|\lambdaeps(r)|\leq\revis{|\lambda(r)|+c(1+|r|)}\leq c(1+|r|)$
\revis{(by \eqref{deflambdaeps} and~\eqref{hpPi})}
and $V\subset\Lx6$, it follows that
\Beq
  \norma{\lambdaeps(\phieps)}_{\L2{\Lx6}} \leq c
  \aand
  \norma{\thetaeps\lambdaeps(\phieps)}_{\L2{\Lx3}} \leq c \,.
  \label{stimelambda}
\Eeq
Therefore, by applying~\eqref{stimaxipuntuale},
we deduce that
\Beq
  \norma\xieps_{\L2H}
  \leq c \, \bigl( 1 + \norma\phieps_{\H1H} + \norma{\thetaeps\lambdaeps(\phieps)}_{\L2H} \bigr)
  \leq c \,.
  \label{stimaxieps}
\Eeq
By comparison in \eqref{secondaeps}, we derive an estimate of $\Delta\phieps$ in~$\L2H$, whence
\Beq
  \norma\phieps_{\L2W} \leq c
  \label{stimaphi}
\Eeq
by \pier{\eqref{bceps} and} elliptic regularity.

\step
Second a priori estimate

The estimate we need \pier{next} is the following
\Bsist
  && \norma\thetaeps_{\H1H\cap\L\infty V}
  + \norma\thetaepsG_{\H1\HG}
  \non
  \\
  && \quad {}
  + \norma\phieps_{\W{1,\infty}H \cap \H1V}
  \leq c \,.
  \label{secondastima}
\Esist
A rigorous proof is given in the Appendix.
Here, we proceed formally.
We take $v=\dt\thetaeps$ and $\vG=\dt\thetaeps\suG$ as test functions in \eqref{primaeps}
and integrate over~$(0,t)$.
At the same time,
we multiply the equation obtained by differentiating \eqref{secondaeps} with respect to time
by $\dt\phieps$ and integrate over~$Q_t$.
Then, we add the equalities just \pier{derived} to each other.
Since the terms involving the product $\dt\thetaeps\,\dt\phieps$ cancel out, we have
\Bsist
  && \intQt |\dt\thetaeps|^2
  + \frac 12 \iO |\nabla\thetaeps(t)|^2
  + \tau \intSt |\dt\thetaepsG|^2
  + \frac \alpha 2 \iG |\thetaepsG(t)|^2
  \non
  \\
  && \quad {}
  + \frac 12 \iO |\dt\phieps(t)|^2
  + \sigma \intQt |\nabla\dt\phieps|^2
  + \intQt \beta'_\eps(\phieps) |\dt\phieps|^2
  \non
  \\
  && = \frac 12 \iO |\nabla\thetaz|^2
  + \frac \alpha 2 \iG |\thetaz\suG|^2
  + \frac 12 \iO |\dt\phieps(0)|^2
  \revis{{} + \intQt \thetaeps \lambdaeps'(\phieps) |\dt\phieps|^2 }
  \non
  \\
  && \quad {}
  + \alpha \intSt mu \, \dt\thetaepsG
  - \intQt \pi'(\phieps) \, |\dt\phieps|^2 .
  \non
\Esist
All the terms on the \lhs\ are nonnegative, \pier{in particular for} 
the last one \pier{we use the monotonicity of $\betaeps$}.
Moreover, the first two terms on the \rhs\ are finite by \eqref{hpdati}
and the last two integrals can be treated in an obvious way
by also taking the \Lip\ continuity of $\pi$ into account.
\revis{In order to control the term involving $\lambdaeps'$, 
we observe that \eqref{deflambdaeps} easily implies that
$|\lambdaeps'(r)|\leq c$ for all $r\in\erre$ and $\eps\in(0,1)$.
We also recall that $V\subset\Lx4$ and term $C$ the norm of the embedding.
Therefore, by owing to the \Holder\ and Young inequalities, we obtain
\Bsist
  && \intQt \thetaeps \lambdaeps'(\phieps) |\dt\phieps|^2
  \leq c \iot \norma{\thetaeps(s)}_4 \, \norma{\dt\phieps(s)}_4\, \norma{\dt\phieps(s)}_2 \, ds
  \non
  \\
  && \leq \frac \sigma {2C^2} \iot \norma{\dt\phieps(s)}_4^2 \, ds
  + c \iot \norma{\thetaeps(s)}_4^2 \, \norma{\dt\phieps(s)}_2^2 \, ds
  \non
  \\
  && \leq \frac \sigma 2 \iot \normaH{\nabla\dt\phieps(s)}^2 \, ds
  + \frac \sigma 2 \iot \normaH{\dt\phieps(s)}^2 \, ds
  + c \iot \normaV{\thetaeps(s)}^2 \, \norma{\dt\phieps(s)}_2^2 \, ds 
  \non
\Esist
and we remark at once that $s\mapsto\normaV{\thetaeps(s)}^2$ is bounded in $L^1(0,T)$ by~\eqref{primastima}.
Finally, we estimate $\dt\phieps(0)$ in~$H$.}
We formally have from~\eqref{secondaeps}
\Beq
  \dt\phieps(0)
  = \sigma \Delta\phiz - \betaeps(\phiz) - \pi(\phiz) + \thetaz \lambdaeps(\phiz) \quad \pier{\aeO}
  \non
\Eeq
and we can owe to \eqref{hpdati}, \eqref{propYosida}, the inequality
$|\lambdaeps(r)|\leq\revis{|\lambda(r)|+c(1+|r|)}\leq c(1+|r|)$ for every $r\in\erre$,
and account for the continuous embedding $V\subset\Lx4$ \revis{once more}.
Then, \eqref{secondastima} follows by the Gronwall lemma.

\step
Third a priori estimate

The \Holder\ inequality, the continuous embedding $V\subset\Lx4$ and \eqref{secondastima} imply
\Beq
  \norma{\thetaeps\lambdaeps(\phieps)}_{\L\infty H}
  \leq c \, \norma\thetaeps_{\L\infty{\Lx4}} \bigl( 1 + \norma\phieps_{\L\infty{\Lx4}} \bigr)
  \leq c .
  \label{nuovalambdaeps}
\Eeq
On the other hand, \eqref{stimaxipuntuale} holds for the approximating solution.
Thus, $\xieps$~is bounded in $\L\infty H$
and a bound for $\Delta\phieps$ in $\L\infty H$ follows by comparison in~\eqref{secondaeps}.
Hence, by also using \pier{\eqref{bceps} and the elliptic regularity theory}, we~have
\Beq
  \norma\phieps_{\L\infty W}
  + \norma\xieps_{\L\infty H}
  \leq c \,.
  \label{terzastima}
\Eeq

\step
Conclusion of the proof

By standard weak and weak star compactness results,
we have for a subsequence
\Bsist
  & \thetaeps \to \theta
  & \quad \hbox{weakly star in $\H1H\cap\L\infty V$}
  \label{convtheta}
  \\
  & \thetaepsG \to \thetaG
  & \quad \hbox{weakly in $\H1\HG$}
  \label{convthetaG}
  \\
  & \phieps \to \phi
  & \quad \hbox{weakly star in $\W{1,\infty}H\cap\H1V\cap\L\infty W$}
  \qquad
  \label{convphi}
  \\
  & \xieps \to \xi
  & \quad \hbox{weakly in $\L\infty H$}
  \label{convxi}
\Esist
and $\thetaG(t)=\theta(t)\suG$ \aat.
By owing to the compact emedding $W\subset\Cx0$ and to \cite[Sect.~8, Cor.~4]{Simon},
we can also assume that the selected subsequence satisfies
\Beq
  \phieps \to \phi
  \quad \hbox{uniformly in $Q$}
  \label{strongphi}
\Eeq
so that $\lambdaeps(\phieps)$ converges to $\lambda(\phi)$ uniformly in $Q$
since $\lambdaeps(r)$ converges to $\lambda(r)$
uniformly on every bounded interval.
We deduce that
$\pier{\lambdaeps}(\phieps)\dt\phieps$ and $\thetaeps\lambdaeps(\phieps)$ converge to
$\lambda(\phi)\dt\phi$ and to $\theta\lambda(\phi)$
at least weakly in~$\LQ2$.
Finally, we have $\xi\in\beta(\phi)$ \aeQ\ by applying, e.g., \cite[\pier{Prop.~2.2, p.~38}]{Barbu}.
Therefore, we can pass to the limit in the integrated version of problem \Pbleps\
written with time dependent test functions
and easily conclude that $(\theta,\thetaG,\phi,\xi)$ is a solution to problem \Pbl.\qed

\medskip

\gianni{
Now, we prove Theorem~\ref{Regularity}. 
For the claim~$i)$, we first consider the Cauchy problem for the linear variational equation
\Bsist
  && \iO \dt\theta \, v
  + \iO \nabla\theta \cdot \nabla v
  + \tau \iG \dt\thetaG \, \vG
  + \alpha \iG \thetaG \, \vG
  = \iO \psi v + \iG \psiG \vG
  \qquad
  \non
  \\
  && \quad \hbox{for every $(v,\vG)\in V\times\VG$ such that $\vG=v\suG$ and \aet}
  \label{primabis}
\Esist
where $\psi$ and $\psiG$ are given, and prove the following}

\Bprop
\label{Thetabdd}
Assume that $\psi\in\L\infty H$, $\psiG\in\pier{\LS\infty}$ and $\thetaz\in V\cap\Lx\infty$
and that the corresponding norms are estimated by some constant~$C$.
Moreover, assume that
$(\theta,\thetaG)$ enjoys the regularity \accorpa{regtheta}{regthetaG}
and satisfies the variational equation~\eqref{primabis} and the initial condition $\theta(0)=\thetaz$.
Then, $\theta$~is bounded and the estimate
\Beq
  \norma\theta_\infty
  \leq \hat C
  \non
\Eeq
holds true with a constant $\hat C$ that depends only on
$\Omega$, $T$, $\tau$, $\alpha$ and~$C$.
\Eprop

\Bdim
This is a regularity result for a linear problem.
Hence, it can be established by considering the following cases:
$$
  a) \quad \psi = 0 \aand \thetaz = 0 \,; \qquad
  b) \quad \psiG = 0 \,.
$$
Let us consider the first one.
For any integer $n>0$ and $p\in(2,+\infty)$ we define
\Beq
  T_n(r) := n \tanh \frac rn
  \aand
  \tilde T_{np}(r) := \int_0^{|r|} \bigl( T_n(s) \bigr)^{p-1} ds
  \quad \hbox{for $r\in\erre$} .
  \non
\Eeq
Now, we set $\thetan:=T_n(\theta)$ and $\thetanG:=T_n(\thetaG)$,
and we test \eqref{primabis} by
$v=|\thetan|^{p-1}\sign\theta$ and $\vG=|\thetanG|^{p-1}\sign\thetaG$,
where the sign function is extended by $\sign(0)=0$.
Then, we integrate over $(0,T)$.
As $\sign\thetanG=\sign\thetaG$, we obtain
\Bsist
  && \iO \tilde T_{np}(\theta(T))
  + \tau \revis\iG \tilde T_{np}(\thetaG(T))
  + (p-1) \intQ |\thetan|^{p-2} \, T'_n(\thetan) |\nabla\theta|^2
  + \alpha \intS |\thetaG| |\thetanG|^{p-1}
  \non
  \\
  && = \intS \psiG \, |\thetanG|^{p-1} \sign\thetaG \,.
  \non
\Esist
All the terms on the \lhs\ are nonnegative.
Hence, by observing that $|r|\geq|T_n(r)|$ for every~$r$,
and applying the Young inequality, we have
\Bsist
  && \alpha \intS |\thetanG|^p
  \leq \alpha \intS |\thetaG| |\thetanG|^{p-1}
  \leq \intS |\psiG| \, |\thetanG|^{p-1}
  \non
  \\
  && = \intS \alpha^{-1/p'}|\psiG| \, \alpha^{1/p'} |\thetanG|^{p-1}
  \leq \intS \Bigl( \frac 1p \, \alpha^{-p/p'} |\psiG|^p + \frac 1{p'} \, \alpha |\thetanG|^p \Bigr)
  \non
\Esist
whence immediately
\Beq
  \frac \alpha p \intS |\thetanG|^p
  \leq \frac 1p \, \alpha^{-p/p'} \intS |\psiG|^p
  \qquad \hbox{or} \qquad
  \alpha^{1/p} \norma{\thetanG}_p
  \leq \alpha^{-1/p'} \norma\psiG_p \,.
  \non
\Eeq
By noting that $|T_n(r)|\Neto|r|$ for every $r\in\erre$ as $n\neto\infty$,
we can let first $n$ and then $p$ tend to infinity and deduce that $\thetaG$ is bounded.
Namely
\Beq
  \norma{\thetaG}_\infty
  \leq \alpha^{-1} \norma\psiG_\infty \leq C /\alpha \,.
\Eeq
Hence, we can apply \cite[Thm.~7.1]{LSU} with $q=2$ and $r=\infty$
and conclude.

In case~$b)$ we adapt the proof of \cite[Thm.~7.1]{LSU}
(still with $q=2$ and $r=\infty$) to~our situation.
For $k\geq\max\{1,C\}$ we set
$\theta^k:=(\theta-k)^+$ and $\thetaG^k:=(\thetaG-k)^+$
and \pier{take $v=\theta^k$ in \eqref{primabis}.}
By simply writing $\thetaG=(\thetaG-k)+k$
and observing that $\theta^k(0)=0$ since $\thetaz\leq C$, we obtain
\Beq
  \frac 12 \iO |\theta^k(t)|^2
  + \frac \tau 2 \iG |\thetaG^k(t)|^2
  + \intQt |\nabla\theta^k|^2
  + \alpha \intSt |\thetaG^k|^2
  + k \alpha \intSt \thetaG^k
  = \intQt \psi\theta^k
  \non
\Eeq
whence also
\Beq
  \frac 12 \iO |\theta^k(t)|^2
  + \intQt |\nabla\theta^k|^2
  \leq \intQt \psi\theta^k .
  \non
\Eeq
This corresponds to formula \cite[(7.6), p.~183]{LSU}
and the whole argument can be performed in the same way
till formula \cite[(7.14), p.~186]{LSU}.
Then, just one modification is needed.
Namely, we account for \cite[Rem.~6.2, p.~103]{LSU}
in applying \cite[Thm.~6.1, p.~102]{LSU}
since no upper bound for $\thetaG$ is known now.
This leads to the desired estimate from above $\theta\leq\hat C$.
The corresponding estimate from below is obtained
by applying the former to~$-\theta$.
\Edim

\gianni{Now, we apply the above result by observing that
the pair $(\theta,\thetaG)$ we are interested in 
satisfies \eqref{primabis} with $\psi:=\lambda(\phi) \dt\phi$ and $\psiG:=mu$,
and notice that these functions
belong to $\L\infty H$ and to $\LS\infty$, respectively.
Moreover, the corresponding norm of $\psi$ has been already estimated by~\eqref{stimasoluz},
while an upper bound of the norm of $\psiG$ is supposed to be given in the statement.
This yields the claim~$i)$ of Theorem~\ref{Regularity}}.

The next step should be the proof of $ii)$ of Theorem~\ref{Regularity}.
To this aim, we just refer to \cite[Thm.~2.2, \gianni{$iii)$}]{CoGiMaRo}.
In fact, the proof given there
shows that the component $\xi$ of a pair $(\phi,\xi)$ satisfying \accorpa{regphi}{regxi}
and solving the homogeneous Neumann problem for the equations
\Beq
  \dt\phi - \sigma \Delta\phi + \xi + \pi(\phi) = f
  \aand
  \xi \in \beta(\phi) \quad \aeQ
  \non
\Eeq
is bounded whenever $f\in\LQ\infty$.
Such a statement is proved just with $f=\theta$ in the quoted paper
by knowing that $\theta$ is bounded.
However, the same proof is valid with any bounded~$f$.
Here, we have $f=\theta\lambda(\phi)$,
and both $\theta$ and $\lambda(\phi)$ are bounded.


\section{Existence of an optimal control}
\label{OPTIMUM}
\setcounter{equation}{0}

We prove Theorem~\ref{Optimum} by the direct method.
Since $\Uad$ is nonempty,
we can take a minimizing sequence $\graffe{\un}$ for the optimization problem
and, for any~$n$, we can consider the corresponding solution $(\thetan,\thetanG,\phin,\xin)$ to problem~\Pbl.
Then, $\graffe{\un}$ is bounded in $\Lx\infty$ and estimate \eqref{stimasoluz}
holds for $(\thetan,\thetanG,\phin,\xin)$.
Therefore, we have
\Bsist
  & \un \to \uopt
  & \quad \hbox{weakly star in $\Lx\infty$}
  \non
  \\
  & \thetan \to \thetaopt
  & \quad \hbox{weakly star in $\H1H \cap \L\infty V$}
  \non
  \\
  & \thetanG \to \thetaoptG
  & \quad \hbox{weakly star in $\H1\HG$}
  \non
  \\
  & \phin \to \phiopt
  & \quad \hbox{weakly star in $\W{1,\infty}H \cap \H1V \cap \L\infty W$}
  \non
  \\
  & \xin \to \xiopt
  & \quad \hbox{weakly star in $\L\infty H$}
  \non
\Esist
at least for a subsequence.
Then, $\uopt\in\Uad$ since $\Uad$ is closed in~$\calX$.
Moreover, the initial conditions for $\thetaopt$ and $\phiopt$ are satisfied
and $\thetaoptG=\thetaopt\suG$ \aet.
Thus, we can easily conclude by standard argument.
Indeed, $\graffe{\phin}$~converges to $\phiopt$ uniformly in~$Q$
due to the compact embedding $W\subset\Cx0$
(see, e.g., \cite[Sect.~8, Cor.~4]{Simon}),
whence $\pi(\phin)$ and $\lambda(\phin)$ converge to $\pi(\phiopt)$ and $\lambda(\phiopt)$
in the same topology.
We also deduce that
$\lambda(\thetan)\dt\phin$ and $\thetan\lambda(\phin)$
converge to $\lambda(\thetaopt)\pier{\dt\phiopt}$ and $\pier{\thetaopt}\lambda(\phiopt)$ at least weakly in~$\LQ2$,
and that
$\xiopt\in\beta(\phiopt)$ (\pier{note that 
$\phin\to \phiopt$ strongly in $\LQ2$ and} 
see, e.g., \cite[\gabri{\gianni{Prop}.~2.1, p.~29}]{Barbu}).
Hence, $(\thetaopt,\thetaoptG,\phiopt,\xiopt)$ satisfies the variational formulation in the integral form of problem~\Pbl\
corresponding to~$\uopt$.
Therefore
\Beq
  \calJ(\uopt) = \calI(\thetaopt,\thetaoptG,\phiopt)
  \leq \liminf_{n\to\infty} \calI(\thetan,\thetanG,\phin)
  = \lim_{n\to\infty} \calJ(\un)
  = \inf_{u\in\Uad}\calJ(u) .
  \non
\Eeq


\section{The control-to-state mapping}
\label{FRECHET}
\setcounter{equation}{0}

As \pier{sketched} in Sections~\ref{STATEMENT},
the main point is the \Frechet\ differentiability of the control-to-state mapping~$\calS$.
This involves the \lineariz ed problem \Linpbl\
and we first prove Proposition~\ref{Existlin},
i.e., well-posedn\gabri{e}ss for the \lineariz ed problem
and the continuous dependence of its solution on the parameter~$h$.
\gianni{It is understood that all the assumptions of Theorem~\ref{Regularity}
as well as \pier{\HPstrutturater} are in force.}

\step
Well-pose\gabri{d}ness

We aim to apply a contraction argument.
To this end, we observe that all the coefficient\gabri{s} that enter the equations but $\overline c:=\lambda'(\phibar)\dt\phibar$
are bounded thanks to Corollary~\ref{Bddaway}
and that the possibly unbounded coefficient $\overline c$ belongs to~$\L\infty H$.
We define the maps $\calF_1$, $\calF_2$ and $\calF$ in a proper functional framework as follows.
For $\Thetabar\in\LQ2$, we consider the problem for $\Phi$
given by \accorpa{linseconda}{linbc} and the second \pier{condition in} \eqref{lincauchy},
where $\Theta$ in \eqref{linseconda} is replaced by~$\Thetabar$.
We obtain a linear parabolic problem
which has a unique solution $\Phi$ satisfying~\eqref{regPhi}.
We set $\calF_1(\Thetabar):=\Phi$.
By doing that, we obtain a map $\calF_1:\LQ2\to\H1H$.
Now, we fix $\Phi\in\H1H$
and consider the problem for $(\Theta,\ThetaG)$ given by
the variational equation \eqref{linprima} and the initial condition $\Theta(0)=0$.
Such a problem has a unique solution $(\Theta,\ThetaG)$ satisfying \accorpa{regTheta}{regThetaG},
as one can see by arguing as in the proof of Lemma~\ref{Phitheta},
and we set $\calF_2(\pier{\Phi}):=\Theta$.
In this way, we obtain a map $\calF_2:\H1H\to\LQ2$.
We set $\calF:=\calF_2\circ\calF_1$ and prove that
some iterated $\calF^m$ of $\calF$ is a contraction in~$\LQ2$
by deriving some estimates involving $\calF_1$ and~$\calF_2$, separately.
For $\Thetabar\in\LQ2$ we write \eqref{linseconda} with $\Thetabar$ in place of $\Theta$
and multiply it by~$\Phi$.
Then, we integrate over $Q_t$ and obtain
\Beq
  \frac 12 \iO |\Phi(t)|^2
  + \sigma \intQt |\nabla\Phi|^2
  = \intQt \bigl( \thetabar \lambda'(\phibar) - \gamma'(\phibar) \bigr) |\Phi|^2
  + \intQt |\lambda(\phibar)| \, |\Thetabar| \, |\Phi| .
  \non
\Eeq
\pier{Hence, with the help of the Young inequality and the Gronwall lemma, we} deduce the a priori estimate
\Beq
  \norma\Phi_{\Lt\infty H\cap\Lt2V}
  \leq C \, \norma\Thetabar_{L^2(Q_t)}
  \quad \hbox{for every $t\in[0,T]$}
  \non
\Eeq
where the constant $C$ we have marked with a capital letter for a future reference
does depend \gabri{neither} on~$t$ nor on~$\Thetabar$.
As $\calF_1$ is linear, this means that
\Beq
  \norma{\calF_1(\Thetabar_1) - \revis{\calF_1}(\Thetabar_2)}_{\Lt\infty H\cap\Lt2V}
  \leq C \, \norma{\Thetabar_1 - \Thetabar_2}_{L^2(Q_t)}
  \quad \hbox{for $t\in[0,T]$}
  \label{stimaFuno}
\Eeq
for every $\Thetabar_1,\Thetabar_2\in\LQ2$.
Now, for given $\Phi_i\in\H1H$, $i=1,2$,
we consider the problems corresponding to the definition of $\Theta_i=\calF_2(\Phi_i)$
and take the difference.
By setting $\Phi:=\Phi_1-\Phi_2$
and analogously defining $\Theta$ and $\ThetaG$ for brevity,
we see that \eqref{linprima} holds with $h$ replaced by~$0$.
We integrate such an equality with respect to time.
By observing that
$\lambda(\phibar)\dt\Phi+\lambda'(\phibar)\dt\phibar\,\Phi=\dt(\lambda(\phibar)\,\Phi)$,
we simply~obtain
\Beq
  \iO \Theta \, v
  + \iO \nabla(1*\Theta) \cdot \nabla v
  + \tau \iG \ThetaG \, \vG
  + \alpha \iG (1*\ThetaG) \, \vG
  = - \iO  \lambda(\phibar) \, \Phi \, v
  \label{intlinprima}
\Eeq
for every $(v,\vG)\in V\times\VG$ such that $\vG=v\suG$ and \aet.
By choosing $v=\Theta$ and $\vG=\ThetaG$,
\gianni{integrating with respect to time}
and forgetting some nonnegative terms on the \lhs, we easily infer~that
\Beq
  \norma{\Theta_1 - \Theta_2}_{\Lt2H}
  \leq D \, \norma{\Phi_1 - \Phi_2}_{\Lt2H}
  \quad \hbox{for $t\in[0,T]$}
  \label{stimaFdue}
\Eeq
where we have marked the constant $D$ for convenience.
If we combine \eqref{stimaFdue} written for $\Phi_i=\calF_1(\Thetabar_i)$
with \eqref{stimaFuno},
we deduce the estimate
\Bsist
  && \norma{\calF(\Thetabar_1) - \calF(\Thetabar_2)}_{\Lt2H}^2
  = \iot \normaH{(\Theta_1 - \Theta_2)(s)}^2 \, ds
  \non
  \\
  && \leq D^2 \iot \normaH{(\Phi_1 - \Phi_2)(s)}^2 \, ds
  \leq D^2 \iot \norma{\gianni{\Phi_1 - \Phi_2}}_{\Ls\infty H}^2 \, ds
  \non
  \\
  && \leq D^2 C^2 \iot \norma{\Thetabar_1 - \Thetabar_2}_{\Ls2H}^2 \, ds
  \quad \hbox{for every $t\in[0,T]$}
  \non
\Esist
and this can be iterated.
By doing that, we obtain the inequality
\Beq
  \norma{\calF^m(\Thetabar_1) - \calF^m(\Thetabar_2)}_{\Lt2H}^2
  \leq \frac {(C^2 D^2 \, t)^m} {m!} \, \norma{\Thetabar_1 - \Thetabar_2}_{\Lt2H}^2
  \non
\Eeq
for every $\Theta_i\in\LQ2$, $t\in[0,T]$ and $m\geq1$.
By choosing $m$ such that $(C^2 D^2 \, T)^m<m!$ and $t=T$,
we see that $\calF^m$ is a contraction in $\LQ2$,
whence $\calF$ has a unique fixed point~$\Theta$.
Then, $\Theta$~and the associated functions $\ThetaG$ and $\Phi$
that enter the construction of $\calF_1$ and $\calF_2$
provide a unique solution to the \lineariz ed problem~\pier{\Linpbl\
with the regularity}~\Reglin.

\step
Continuous dependence

This is given by \pier{the estimate~\eqref{stimaFrechet} we} prove at once.
By integrating \eqref{linprima} with respect to time
and proceeding as we did for~\eqref{intlinprima}, we have
\Beq
  \iO \Theta \, v
  + \iO \nabla(1*\Theta) \cdot \nabla v
  + \tau \iG \ThetaG \, \vG
  + \alpha \iG (1*\ThetaG) \, \vG
  = - \iO  \lambda(\phibar) \, \Phi \, v
  + \alpha \iG m (1*h) \, \vG
  \non
\Eeq
for every $(v,\vG)\in V\times\VG$ such that $\vG=v\suG$ and \aet.
Now, we \pier{take $v=\Theta$, $\vG=\ThetaG$}
and integrate over~$(0,t)$.
Besides, we multiply \eqref{linseconda} by $\delta\dt\Phi$,
where $\delta$ is a positive parameter, and integrate over~$Q_t$.
Then, we add the equalities we obtain to each other.
\revis{Furthermore, we add the same term $(\delta/(4T))\iO|\Phi(t)|^2$ to both sides for convenience.}
By also owing to the boundedness of the coefficients and to the Young inequality,
\revis{and setting $C:=\norma{\lambda(\phibar)}_\infty^2\,$,}
we~have
\Bsist
  && \intQt |\Theta|^2
  + \frac 12 \iO |\nabla(1*\Theta)(t)|^2
  + \tau \intSt |\ThetaG|^2
  + \frac \alpha 2 \iG |(1*\ThetaG)(t)|^2
  \non
  \\
  && \quad {}
  + \delta \intQt |\dt\Phi|^2
  + \frac {\delta\sigma} 2 \iO |\nabla\Phi(t)|^2
  \gianni{{}+ \frac \delta {4T} \iO |\Phi(t)|^2}
  \non
  \\
  \separa
  && = \revis{- \intQt \lambda(\phibar) \Phi \Theta}
  + \alpha \intSt m (1*h) \, \ThetaG
  \non
  \\
  && \quad {}
  + \delta \intQt \bigl( - \gamma'(\phibar) + \thetabar \pier{\lambda'} (\phibar) \bigr) \Phi \, \dt\Phi
  + \delta \intQt \lambda(\phibar) \, \Theta \, \dt\Phi 
  \revis{{}+ \frac \delta {4T} \iO |\Phi(t)|^2}
  \non
  \\
  && \leq \revis{\delta \, C \intQt |\Theta|^2 + \frac 1 {4\delta} \intQt |\Phi|^2}
  \revis{{} + {}} \frac \tau 2 \intSt |\ThetaG|^2
  + c \intSt |1*h|^2
  \non
  \\
  && \quad {}
  + \frac \delta 2 \intQt |\dt\Phi|^2
  + \delta \, c \intQt |\Phi|^2
  + \delta \, C \intQt |\Theta|^2
  + \frac \delta {4T} \iO |\Phi(t)|^2 .
  \non
\Esist
\revis{Moreover, by} arguing as we did to derive the first \pier{inequality in}~\eqref{auxil},
we see~that
\Beq
  \frac \delta {4T} \iO |\Phi(t)|^2
  \leq \frac \delta 4 \intQt |\dt\Phi|^2 .
  \non
\Eeq
Therefore, by choosing $\delta$ such that $\revis 2\delta C<1$ and applying the Gronwall lemma, 
we~obtain
\Beq
  \norma\Theta_{\LQ2}
  + \norma{1*\Theta}_{\L\infty V}\
  + \norma\ThetaG_{\LS2}
  + \norma\Phi_{\H1H\cap\L\infty V}
  \leq c \norma{1*h}_{\LS2} \,.
  \non
\Eeq
\pier{At this point, it is easy to see that the above estimate} implies~\eqref{stimaFrechet}.\qed

Here is the main result of this section.

\Bthm
\label{Fdiff}
Let $\ubar\in\calX$.
Then, $\calS$~is \Frechet\ differentiable at~$\ubar$
and the \Frechet\ derivative $[D\calS](\ubar)$
\pier{is precisely} the map $\calD\in\calL(\calX,\calY)$ defined
in the statement of Proposition~\ref{Existlin}.
\Ethm

\Bdim
We fix $\ubar\in\calX$ and the corresponding state $(\thetabar,\thetabarG,\phibar):=\calS(\ubar)$
and, for $h\in\calX$, we~set
\Beq
  (\thetah,\thetahG,\phih) := \calS(\ubar+h)
  \aand
  (\zetah,\zetahG,\etah) := (\thetah-\thetabar-\Theta,\thetahG-\thetabarG-\ThetaG,\phih-\phibar-\Phi)
  \non
\Eeq
where $(\Theta,\ThetaG,\Phi)$ is the solution to the linearized problem corresponding to~$h$.
We have to prove that
$\norma{(\zetah,\zetahG,\etah)}_\calY/\norma h_\calX$ tends to zero as $\norma h_\calX$ tends to zero.
More precisely, we show that some constant $c$ exists such~that
\Beq
  \norma{(\zetah,\zetahG,\etah)}_\calY
  \leq c \norma h_{\revis{\LS2}}^2
  \quad \hbox{provided that $\norma h_\calX\leq1$}
  \label{tesiFrechet}
\Eeq
and this is even stronger than necessary.
So, we assume $\norma h_\infty\leq 1$
and make a preliminary observation.
As $\ubar$ is fixed and $\norma{\ubar+h}_\infty\leq\norma\ubar_\infty+1$
for every $h$ under consideration,
we can apply Theorem~\ref{Regularity} and Corollary~\ref{Bddaway}
and find constants $\thetamin,\thetamax\in\erre$ and
$\phimin,\phimax\in D(\beta)$ independent of $h$ such that
\Bsist
  && \thetamin \leq \thetabar \leq \thetamax
  \aand
  \thetamin \leq \thetah \leq \thetamax
  \quad \aeQ
  \label{perTaylortheta}
  \\
  && \phimin \leq \phibar \leq \phimax
  \aand
  \phimin \leq \phih \leq \phimax
  \quad \aeQ .
  \label{perTaylorphi}
\Esist
Moreover, we can also \pier{exploit} the second part of Theorem~\ref{Wellposedness}
with $u_1=\ubar+h$ and $u_2=\ubar$.
By doing that and also owing to the continuous embedding $V\subset\Lx4$,
we derive from~\eqref{contdep}
\Bsist
  && \gianni{\norma{\thetah-\thetabar}_{\L\infty H}} \leq C' \norma h_{\LS2} \,, \quad
  \norma{\phih-\phibar}_{\LQ4}
  \leq MC' \norma h_{\LS2}
  \label{contdepA}
  \\
  && \norma{\phih-\phibar}_{\L\infty{\Lx4}}
  \leq MC' \norma h_{\LS2}
  \label{contdepB}
\Esist
\Accorpa\Dacontdep contdepA contdepB
where the constants $M$ and $C'$ do not depend on~$h$.
Now, the problem solved by~$(\zetah,\zetahG,\etah)$ is the following
\Bsist
  && \iO \dt\zetah \, v
  + \iO \nabla\zetah \cdot \nabla v
  + \tau \iG \dt\zetahG \, \vG
  + \alpha \iG \zetahG \vG
  \non
  \\
  && = - \iO \bigl\{
    \lambda(\phih) \dt\phih
    - \lambda(\phibar) \dt\phibar
    - \lambda(\phibar) \dt\Phi
    - \lambda'(\phibar) \dt\phibar \, \revis\Phi
  \bigr\} v
  \non
  \\
  && \quad \hbox{for every $(v,\vG)\in V\times\VG$ such that $\vG=v\suG$ and \aet}
  \label{primah}
  \\
  \separa
  && \dt\etah - \sigma \Delta\etah
  = - E_1^h + E_2^h
  \qquad \aeQ, \ \ \hbox{where}
  \label{secondah}
  \\
  && \quad E_1^h :=
    \gamma(\phih)
    - \gamma(\phibar)
    - \gamma'(\phibar) \, \Phi
  \label{defEunoh}
  \\
  && \quad E_2^h :=
    \thetah \lambda(\phih)
    - \thetabar \lambda(\phibar)
    - \lambda(\phibar) \Theta
    - \thetabar \lambda'(\phibar) \, \Phi
  \label{defEdueh}
  \\
  \separa
  && \dn\etah = 0 \quad \aeS
  \aand
  \zetah(0) = \etah(0) = 0 \quad \aeO.
  \label{bcich}
\Esist
By integrating \eqref{primah} with respect to time, we obtain
\Bsist
  && \iO \zetah \, v
  + \iO \nabla(1*\zetah) \cdot \nabla v
  + \tau \iG \zetahG \, \vG
  + \alpha \iG (1*\zetahG) \vG
  = - \iO Z^h v
  \qquad
  \non
  \\
  && \quad \hbox{for every $(v,\vG)\in V\times\VG$ such that $\vG=v\suG$ and \aet}
  \label{intprimah}
  \\
  && \hbox{where} \quad
  Z^h := \Lambda(\phih) - \Lambda(\phibar) - \lambda(\phibar) \Phi \,.
  \label{defZh}
\Esist
At this point, we \pier{take  $v=\zetah$ and $\vG=\zetahG$ in \eqref{intprimah}}
and integrate over~$(0,t)$.
Besides, we multiply \eqref{secondah} by $\etah$ and integrate over~$Q_t$.
Finally, we add the resulting equalities to each other.
We obtain
\Bsist
  && \intQt |\zetah|^2
  + \frac 12 \iO |\nabla(1*\zetah)(t)|^2
  + \tau \intSt |\zetahG|^2
  + \frac \alpha 2 \iG |(1*\zetahG)(t)|^2
  \non
  \\
  && \quad {}
  + \frac 12 \iO |\etah(t)|^2
  + \sigma \intQt |\nabla\etah|^2
  \non
  \\
  \separa
  && = - \intQt Z^h \zetah
  - \intQt E_1^h \etah
  + \intQt E_2^h \etah
  \label{usareTaylor}
\Esist
and we now estimate each term on the \rhs.
To this end, we represent the functions $Z^h$ and $E_i^h$\pier{, $i=1,2$,} in different forms.
By applying the Taylor expansion to $\Lambda$ and~$\gamma$,
we find functions $\hat\phi$ and $\tilde\phi$
taking values between the ones of $\phibar$ and $\phih$ such~that
\Bsist
  && \Lambda(\phih) - \Lambda(\phibar)
  = \lambda(\phibar) \, (\phih - \phibar)
  + \frac 12 \, \lambda'(\hat\phi) \, (\phih-\phibar)^2
  \non
  \\
  && \gamma(\phih) - \gamma(\phibar)
  = \gamma'(\phibar) \, (\phih - \phibar)
  + \frac 12 \, \gamma''(\tilde\phi) \, (\phih-\phibar)^2
  \non
\Esist
whence immediately
\Beq
  Z^h = \lambda(\phibar) \etah
  + \frac 12 \, \lambda'(\hat\phi) \, (\phih-\phibar)^2
  \aand
  E_1^h = \gamma'(\phibar) \etah
  + \frac 12 \, \gamma''(\tilde\phi) \, (\phih-\phibar)^2 \,.
  \non
\Eeq
In order to rewrite $E_2^h$, we observe that,
if $G:\erre^2\to\erre$ is a $C^2$~function,
we have for every $(y,z)$ and $(y',z')$ in~$\erre^2$
and for a suitable $(\breve y,\breve z)$ in between
\Beq
  G(y',z') - G(y,z)
  - \nabla_{yz} G(y,z) \cdot (y'-y,z'-z)
  = \frac 12 \, [y'-y,z'-z] \, D^2_{yz}(\breve y,\breve z) \, ^t[y'-y,z'-z] \,.
  \non
\Eeq
In the particular case $G(y,z)=y\lambda(z)$, the above formula becomes
\begin{align}
 & y' \lambda(z') - y \lambda(z)
  - \lambda(z) (y'-y) - y \lambda'(z) (z'-z)
  \non
  \\
  &= \frac 12 \, \breve y \pier{\lambda''(\breve z)}  (z'-z)^2
  + \lambda'(\breve z) (y'-y) (z'-z) .
  \non
\end{align}
Therefore, there exist functions $\breve\theta$ and $\breve\phi$
taking values between the ones of $\thetabar$ and~$\thetah$
and the ones of $\phibar$ and~$\phih$, respectively,
such~that
\begin{align}
  &\thetah \lambda(\phih) - \thetabar \lambda(\phibar)
  \non
  \\
  &= \lambda(\phibar) (\thetah-\thetabar)
  + \thetabar \lambda'(\phibar) (\phih-\phibar)
  + \frac 12 \, \breve\theta \pier{\lambda''(\breve \phi)}(\phih-\phibar)^2
  + \lambda'(\pier{\breve \phi}) (\thetah-\thetabar) (\phih-\phibar)
  \non
\end{align}
so that
\Beq
  E_2^h
  = \lambda(\phibar) \, \zetah
  + \thetabar \lambda'(\phibar) \, \etah
   + \frac 12 \, \breve\theta \pier{\lambda''(\breve \phi)}(\phih-\phibar)^2
  + \lambda'(\pier{\breve \phi}) (\thetah-\thetabar) (\phih-\phibar) .
  \non
\Eeq
Notice that $\thetamin\leq\breve\theta\leq\thetamax$
and that $\phimin\leq\hat\phi,\tilde\phi,\breve\phi\leq\phimax$.
Therefore, \pier{recalling \eqref{hpnonlinbis},} the \rhs\ of \eqref{usareTaylor} can be estimated
by the Young inequality as follows
\Bsist
  && - \intQt Z^h \zetah
  - \intQt E_1^h \etah
  + \intQt E_2^h \etah
  \leq \frac 12 \intQt |\zetah|^2
  + c \intQt |\etah|^2
  \non
  \\
  && \quad {}
  + c \intQt |\phih-\phibar|^2 |\zetah|
  + c \intQt |\phih-\phibar|^2 |\etah|
  + c \intQt |\thetah-\thetabar| \, |\phih-\phibar| \, |\etah| \,.
  \non
  \\
  && \leq \frac 34 \intQt |\zetah|^2
  + c \intQt |\etah|^2
  + c \intQ |\phih-\phibar|^4
  + c \intQt |\thetah-\thetabar| \, |\phih-\phibar| \, |\etah| \,.
  \non
  \non
\Esist
On the other hand, \pier{thanks to the \Holder\ inequality and to the continuous embedding $V\subset\Lx4$},
we have for every $\delta>0$
\Bsist
  && \intQt |\thetah-\thetabar| \, |\phih-\phibar| \, |\etah|
  \leq \iot \norma{(\thetah-\thetabar)(s)}_2 \, \norma{(\phih-\phibar)(s)}_4 \, \norma{\etah(s)}_4 \, ds
  \non
  \\
  && \leq \delta \iot \normaV{\etah(s)}^2 \, ds
  + c_\delta \, \norma{\phih-\phibar}_{\L\infty{\Lx4}}^2 \, \norma{\thetah-\thetabar}_{\L2H}^2 \,.
  \non
\Esist
At this point, we choose $\delta$ small enough, apply the Gronwall lemma
and account for~\Dacontdep\ in order to estimate the norms of
$\thetah-\thetabar$ and of $\phih-\phibar$.
This yields \eqref{tesiFrechet} and the proof is complete.
\Edim


\section{Necessary optimality conditions}
\label{OPTIMALITY}
\setcounter{equation}{0}

In this section, we derive the optimality condition~\eqref{cnoptadj}
stated in Theorem~\ref{CNoptadj}.
We start from~\eqref{precnopt} and first prove~\eqref{cnopt}.

\Bprop
\label{CNopt}
Let $\uopt$ be an optimal control and $(\thetaopt,\thetaoptG,\phiopt):=\calS(\uopt)$.
Then, condition \eqref{cnopt} holds.
\Eprop

\Bdim
As already said in Section~\ref{STATEMENT},
we just have to apply the chain rule for \Frechet\ derivatives.
Clearly, the \Frechet\ derivative $[D\calI](\thetabar,\thetabarG,\phibar)$
of the functional $\calI$ exists at every point of $\calY$ and it is given~by
\Beq
  [D\calI](\thetabar,\thetabarG,\phibar):
  (h_1,h_2,h_3)\in\calY \mapsto
  \kuno \intQ (\thetabar - \thetaQ) \, h_1
  + \kdue \iO (\phibar(T) - \phiO) \, h_3(T) \,.
  \non
\Eeq
In particular, this holds if
$(\thetabar,\thetabarG,\phibar)=(\thetaopt,\thetaoptG,\phiopt)=\calS(\uopt)$.
Therefore, Theorem~\ref{Fdiff} and the chain rule ensure that
$\calJ$ is \Frechet\ differentiable at $\uopt$
and that its \Frechet\ derivative $[D\calJ](\uopt)$ at any optimal control $\uopt$
acts as follows
\Beq
  [D\calJ](\uopt):
  h \in \calX \mapsto
  \kuno \intQ (\thetaopt - \thetaQ) \, \Theta
  + \kdue \iO (\phiopt(T) - \phiO) \, \Phi(T)
  \non
\Eeq
where $(\Theta,\ThetaG,\Phi)$ is the solution to the linearized problem corresponding to~$h$.
Therefore, \eqref{cnopt} immediately follows from~\eqref{precnopt}.
\Edim

The next step is the proof of Theorem~\ref{Existenceadj}.
For convenience, we consider the equivalent forward problem
in the unknown $(\ptilde,\ptilde_\Gamma,\qtilde)$
given by $(\ptilde,\ptilde_\Gamma,\qtilde)(t):=(p,\pG,q)(T-t)$
and corresponding to the new coefficient and given terms defined accordingly.
However, to simplify notations, we write $p$, $\pG$ and $q$
instead of $\ptilde$, $\ptilde_\Gamma$ and~$\qtilde$ in the sequel.
The new problem is to find $(p,\pG,q)$ satisfying \Regadj\ and solving
\Bsist
  && \iO \dt p \, v
  + \iO \nabla p \cdot \nabla v
  + \iO a q v
  + \tau \iG \dt\pG \, \vG
  + \alpha \iG \pG \vG
  = \iO f v
  \qquad
  \non
  \\
  && \quad \hbox{for every $v\in V$, $\vG:=v\suG$, and \aet}
  \label{primanew}
  \\
  && \iO \dt q \, v
  + \sigma \iO \nabla q \cdot \nabla v
  + \iO b q v
  = \iO g \dt p \, v
  \qquad
  \non
  \\
  && \quad \hbox{for every $v\in V$ and \aet}
  \label{secondanew}
  \\
  && p(0) = 0
  \aand
  q(0) = \qz
  \quad \aeO
  \label{cauchynew}
\Esist
\Accorpa\Pblnew primanew cauchynew
\gianni{where $a$, \pier{$f$, $b$, $g$} and $\qz$
are deduced from \accorpa{primaadj}{cauchyadj}.
Thus, we have
$a,\,b,\,g\in\LQ\infty$, $f\in\LQ2$ and $\pier{\qz\in V}$.}
We aim to use a contraction argument in a weaker functional framework.
However, it will be clear from the proof
that the unique solution we find satisfies \Regadj.

\step
The equivalent \pier{fixed} point problem

For a given $\qbar\in\LQ2$, we consider
the problem obtained by writing \eqref{primanew} with $q$ replaced by $\qbar$
and the initial condition $p(0)=0$.
If we introduce the spaces $\calV$ and $\calH$ and the bilinear form $a$
given by \accorpa{calVH}{defforma},
and argue as in the proof of Lemma~\ref{Phitheta},
we see that the problem under consideration has a unique solution $(p,\pG)$
satisfying
$p\in\H1H\cap\L\infty V$ and~\eqref{regpG}.
However, 
\pier{$p$ is smoother
since one can argue as} in Remark~\ref{Strongsoluz}
to derive $\Delta p\in\LQ2$ and $\dn p\in\LS2$ from the regularity already achieved.
We set $\pier{\calF_1(\qbar):=p}$ and $\tilde\calF_1(\qbar):=(p,\pG)$.
By doing that, we find a map $\calF_1:\LQ2\to\H1H$
and an associated map $\pier{\tilde\calF_1}$ that we use in the rest of the proof.
Now, for $p\in\H1H$, we consider \eqref{secondanew}
complemented by the second \pier{initial condition in}~\eqref{cauchynew}.
\pier{As $b \in\LQ\infty$, $g\dt p \in\LQ2$ and $\qz\in V$,}
such a problem has a unique solution $q$ satisfying~\eqref{regpq},
and we set $q:=\calF_2(p)$.
We thus obtain a map $\calF_2:\H1H\to\LQ2$
and consider $\calF:\LQ2\to\LQ2$ defined by $\calF:=\calF_2\circ\calF_1$.
\pier{Let us point out that, for a given $\qbar\in\LQ2$, $\calF (\qbar)$ actually takes values in 
the space $\H1H \cap \L\infty V \cap \L2W.$}
The problem under consideration is equivalent to
the existence of a unique fixed point for~$\calF$.
Indeed, if~$q$ is such a fixed point,
$q$~and the corresponding $p$ and $\pG$ provide a solution satisfying \Regadj\ by construction.
Conversely, if $(p,\pG,q)$ solves the problem, then $q$ is a fixed point of~$\calF$.

\step
The contraction argument

It suffices to find a constant $C$ such that
\Beq
  \normaH{\bigl(\calF(\qbar_1)-\calF(\qbar_2)\bigr)(t)}
  \leq C \, \norma{\qbar_1-\qbar_2}_{\Lt2H}
  \label{contrazione}
\Eeq
for every $\qbar_1\,,\,\qbar_2\in\LQ2$ and every $t\in[0,T]$.
Indeed, this implies that some iterated $\calF^m$ of $\calF$ is a contraction.
In order to prove~\eqref{contrazione},
we take any pair of functions $\qbar_i\in\LQ2$, $i=1,2$,
consider the pairs $(p_i,p_{i,\Gamma}):=\pier{\tilde\calF_1}(\qbar_i)$
and the functions $q_i:=\calF(\qbar_i)$
and write their definitions, i.e.,
\eqref{primanew} written with $\qbar_i$ in place of $q$
and~\eqref{secondanew} with $p_i$ in place of~$p$.
Then we take the two differences.
We set for convenience $\qbar:=\qbar_1-\qbar_2$ and similarly introduce \pier{$p, \, \pG$} and~$q$.
At this point, we formally test the first difference by $v=\dt p$ and $\vG=\dt\pG$
(by~avoiding the technical approximation of such test functions by $V$- and $\VG$-valued functions)
and~integrate over~$(0,t)$.
At the same time, we multiply the second one by~$q$ and integrate over~$Q_t$.
We~have
\Bsist
  && \intQt |\dt p|^2
  + \frac 12 \iO |\nabla p(t)|^2
  + \tau \revis\intSt |\dt\pG|^2 
  + \frac \alpha 2 \iG |\pG(t)|^2
  = \pier{- \intQt a \, \qbar \, \dt p}
  \non
  \\
  \noalign{\allowbreak\noindent as well as}
  && \frac 12 \iO |q(t)|^2
  + \sigma \intQt |\nabla q|^2
  = - \intQt b \, |q|^2
  + \intQt g \, \dt p \, q \,.
  \non
\Esist
As $a$, $b$ and $g$ are bounded functions,
it is \sfw\ to deduce that
\Bsist
  && \norma p_{\Ht1H\cap\Lt\infty V}
  + \norma\pG_{\Lt2\HG}
  \leq C_1 \, \norma\qbar_{\Lt2H}
  \non
  \\
  \noalign{\smallskip}
  && \norma q_{\Lt\infty H\cap\Lt2V}
  \leq C_2 \, \norma p_{\Ht1H}
  \non
\Esist
for every $t\in[0,T]$ and some constants $C_1$ and $C_2$
independent of $\qbar_i$ and~$t$.
By combining such inequalities, we deduce that
\Beq
  \norma q_{\Lt\infty H}
  \leq C_1C_2 \, \norma\qbar_{\Lt2H}
  \non
\Eeq
\pier{whence \eqref{contrazione} follows} with $C:=C_1C_2$.
Thus, Theorem~\ref{Existenceadj} is completely proved.\qed

\medskip

At this point, we are ready to prove Theorem~\ref{CNoptadj} on optimality,
i.e., the necessary condition \eqref{cnoptadj} for
$\uopt$ to be an optimal control in terms of the solution $(p,\pG,q)$
of the adjoint problem~\Pbladj.
So, let $\uopt$ be an optimal control and fix an arbitrary $u\in\Uad$.
We write both the variational formulations of the \lineariz ed problem
at the optimal state $(\thetaopt,\thetaoptG,\phiopt):=\calS(\uopt)$
corresponding to $h=u-\uopt$
and the adjoint system.
We have \aet
\begin{align}
  & \iO \dt\Theta \, v
  + \iO \nabla\Theta \cdot \nabla v
  + \iO \lambda(\phiopt) \dt\Phi \, v
  + \iO \lambda'(\phiopt) \dt\phiopt \, \Phi \, v
  \non
  \\
  & \quad {}
  + \tau \iG \dt\ThetaG \, \vG
  + \alpha \iG \ThetaG \, \vG
  = \alpha \iG m (u-\uopt) \, \vG
  \label{vprima}
  \\
  \separa
  & \iO \dt\Phi \, v
  + \sigma \iO \nabla\Phi \cdot \nabla v
  + \iO \gamma'(\phiopt) \, \Phi v
  = \iO \thetaopt \lambda'(\phiopt) \, \Phi v
  + \iO \lambda(\phiopt) \, \Theta v
  \label{vseconda}
  \\
  \separa
  & - \iO \dt p \, v
  + \iO \nabla p \cdot \nabla v
  - \iO \lambda(\phiopt) q v
  - \tau \iG \dt\pG \, \vG
  + \alpha \iG \pG \vG
  \non
  \\
  & \quad{}= \kuno \iO (\thetaopt-\thetaQ) v
  \label{vprimaadj}
  \\
  \separa
  & - \iO \dt q \, v
  + \sigma \iO \nabla q \cdot \nabla v
  + \iO \bigl( \gamma'(\phiopt) - \thetaopt \lambda'(\phiopt) \bigr) q v
  = \iO \lambda(\phiopt) \dt p \, v
  \label{vsecondaadj}
\end{align}
where \eqref{vprima} and \eqref{vprimaadj} hold
for every $(v,\vG)\in V\times\VG$ such that $\vG=v\suG$ ,
while \eqref{vseconda} and \eqref{vsecondaadj} are required for every $v\in V$.
Moreover, the functions at hand satisfy the homogeneous initial or final conditions
as specified in \eqref{lincauchy} and~\eqref{cauchyadj}.
We choose $(v,\vG)=(p,\pG)$ in~\eqref{vprima},
$v=q$ in~\eqref{vseconda},
$(v,\vG)=(-\Theta,-\ThetaG)$ in~\eqref{vprimaadj}
and $v=-\Phi$ in~\eqref{vsecondaadj}.
Then, we sum all the resulting equalities to each other and integrate over~$(0,T)$.
Several terms cancel out and we obtain
\Bsist
  && \intQ \bigl\{
    \dt\Theta \, p
    + \lambda(\phiopt) \dt\Phi \, p
    + \lambda'(\phiopt)\dt\phiopt \, \Phi p
    + \dt\Phi \, q
    + \dt p \, \Theta
    + \dt q \, \Phi + \lambda(\phiopt) \dt p \, \Phi
  \bigr\}
  \non
  \\
  && \quad {}
  + \tau \intS \bigl( \dt\ThetaG \, \pG + \dt\pG \, \ThetaG \bigr)
  \non
  \\
  \separa
  && = \alpha \intS m (u-\uopt) \pG
  - \kuno \intQ (\thetaopt - \thetaQ) \Theta \,.
  \non
\Esist
As the expression between braces is equal to
$\dt\bigl(\Theta p+\lambda(\phiopt)\Phi p+\Phi q\bigr)$,
the above equality becomes
\Bsist
  && \iO \bigl(\Theta p+\lambda(\phiopt)\Phi p+\Phi q\bigr)(T)
  - \iO \bigl(\Theta p+\lambda(\phiopt)\Phi p+\Phi q\bigr)(0)
  \non
  \\
  && \quad {}
  + \tau \intS (\ThetaG\pG)(T)
  - \tau \intS (\ThetaG\pG)(0)
  \non
  \\
  \separa
  && = \alpha \intS m (u-\uopt) \pG
  - \kuno \intQ (\thetaopt - \thetaQ) \Theta \,.
  \non
\Esist
\pier{Owing} to the relations $\ThetaG(0)=\Theta(0)\suG$ and $\pG(T)=p(T)\suG$,
and accounting for the initial conditions \eqref{lincauchy}
and the final conditions~\eqref{cauchyadj},
the above equation reduces~to
\Bsist
  && \iO \Phi(T) \, \kdue \bigl( \phiopt(T) - \phiO \bigr)
  \non
  \\
  && = \alpha \intS m (u-\uopt) \pG
  - \kuno \intQ (\thetaopt - \thetaQ) \Theta \,.
  \non
\Esist
Therefore, the inequality \eqref{cnopt} we have already established in Proposition~\ref{CNopt} implies
\Beq
  \alpha \intS m (u-\uopt) \pG
  \geq 0
  \quad \hbox{for every $u\in\Uad$} \,.
  \non
\Eeq
Moreover, the positive coefficient $\alpha$ can be removed.
At this point, a standard argument leads to the pointwise relation~\eqref{cnoptadj}
and to its consequences listed in the statement,
and the proof of Theorem~\ref{CNoptadj} is complete.


\def\thetah{\theta_\eps^h}
\def\thetahG{\theta_{\eps,\Gamma}^h}
\def\phih{\phi_\eps^h}
\def\xih{\xi_\eps^h}

\section{Appendix}
\label{APPENDIX}
\setcounter{equation}{0}

This section is devoted to a rigorous proof of~\eqref{secondastima}.
With respect to the formal procedure of Section~\ref{STATE},
we replace derivatives with difference quotients, essentially.
For $h\in(0,T)$ and $v\in\LQ2$ or $v\in\LS2$,
we define $v^h$ on $(0,T-h)$ by setting $v^h(t):=v(t+h)$.
We integrate \eqref{primaeps} with respect to time over $(s,s+h)$
and test the equality we obtain by $v=(\thetah-\thetaeps)(s)$ and $\vG=(\thetahG-\thetaepsG)(s)$.
At the same time, we write \eqref{secondaeps} at times $s+h$ and~$s$,
multiply the difference by $(\phih-\phieps)(s)$
and integrate over $\Omega$ with respect to space.
Finally, we add the equalities obtained this way to each other
and integrate over~$(0,t)$ with respect to~$s$.
We~have
\Bsist
  && \intQt |\thetah - \thetaeps|^2
  + \frac 12 \iO |\nabla \bigl( (1*\thetaeps)^h - (1*\thetaeps) \bigr)(t)|^2
  \non
  \\
  && \quad {}
  + \tau \intSt |\thetahG - \thetaepsG|^2
  + \frac \alpha 2 \iG |\bigl( (1*\thetaepsG)^h - (1*\thetaepsG) \bigr)(t)|^2
  \non
  \\
  && \quad {}
  + \frac 12 \iO |(\phih - \phieps)(t)|^2
  + \sigma \intQt |\nabla(\phih - \phieps)|^2
  + \intQt (\xih - \xieps) (\thetah - \thetaeps)
  \non
  \\
  \separa
  && =
  - \intQt \bigl( \pi(\phih) - \pi(\phieps) \bigr) (\phih - \phieps)
  + \intSt m \bigl( (1*u)^h - (1*u) \bigr) (\thetahG - \thetaepsG)
  \non
  \\
  && \quad {}
  + \intQt \bigl\{
    \bigl( \thetah \lambdaeps(\phih) - \thetaeps \lambdaeps(\phieps) \bigr) (\phih - \phieps)
    - \bigl( \Lambdaeps(\phih) - \Lambdaeps(\phieps) \bigr) (\thetah - \thetaeps)
  \bigr\}
  \non
  \\
  && \quad {}
  + \frac 12 \iO |\nabla(1*\thetaeps)(h)|^2
  + \frac \alpha 2 \iG |(1*\thetaepsG)(h)|^2
  + \frac 12 \iO |\phieps(h) - \phiz|^2 .
  \label{persecondastima}
\Esist
All the terms on the \lhs\ are nonnegative, the last one by monotonicity,
and the first integral on the \rhs\ can be \revis{estimated} in an obvious way
by using the \Lip\ continuity of~$\pi$.
Moreover, we~have
\Bsist
  && \intSt m \bigl( (1*u)^h - (1*u) \bigr) (\thetahG - \thetaepsG)
  \non
  \\
  && \leq \frac \tau 2 \intSt |\thetahG - \thetaepsG|^2
  + c \intSt |(1*u)^h - (1*u)|^2
  \non
  \\
  \separa
  && \leq \frac \tau 2 \intSt |\thetahG - \thetaepsG|^2
  + c \, h^2 \, \norma u_{\LS2}^2
  = \frac \tau 2 \intSt |\thetahG - \thetaepsG|^2
  + c \, h^2 \,.
  \non
\Esist
The critical term is the next one,
since the cancellation of the formal procedure
does not occur here.
We observe at once that the functions $\lambdaeps$ have a common \Lip\ constant
since $\lambda$ is \Lip\ continuous.
We apply the mean value theorem to $\Lambdaeps$
and obtain a function $\tilde\phieps$, taking values
between the ones of $\phih$ and $\phieps$, in order~that
\Bsist
  && \intQt \bigl\{
    \bigl( \thetah \lambdaeps(\phih) - \thetaeps \lambdaeps(\phieps) \bigr) (\phih - \phieps)
    - \bigl( \Lambdaeps(\phih) - \Lambdaeps(\phieps) \bigr) (\thetah - \thetaeps)
  \bigr\}
  \non
  \\
  \separa
  && = \intQt \bigl\{
    \bigl( \lambdaeps(\phih) - \lambdaeps(\tilde\phieps) \bigr) (\thetah - \thetaeps) (\phih - \phieps)
    - \thetaeps \bigl( \lambdaeps(\phih) - \lambdaeps(\phieps) \bigr) (\phih - \phieps)
  \bigr\}
  \non
  \\
  && \leq c \intQt |\thetah-\thetaeps| \, |\phih - \phieps|^2
  + c \intQt |\thetaeps| \, |\phih - \phieps|^2 .
  \non
\Esist
We treat the last two integrals, separately, by owing to the \Holder\ and Young inequalities
and to the continuous embedding $V\subset\Lx4$.
As such embedding is also compact, \eqref{primastima} and \cite[Sect.~8, Cor.~4]{Simon} imply that
the functions $\phieps$ are equicontinuous $\Lx4$--valued functions.
Hence, for every $\delta>0$, there exists $h_\delta>0$ such that
$\norma{\phih(s)-\phieps(s)}_4\leq\delta$ for every $\eps\in(0,1)$ and $s\in[0,T-h]$ whenever~$h\leq h_\delta$.
In the sequel, $\delta$~is a positive parameter, say $\delta\in(0,1)$,
whose value is choosen later on,
and it is understood that $h$ does not exceed the corresponding~$h_\delta$.
Therefore, we~have
\Bsist
  && \intQt |\thetah-\thetaeps| \, |\phih - \phieps|^2
  \leq \iot \norma{(\thetah - \thetaeps)(s)}_2 \, \norma{(\phih - \phieps)(s)}_4^2 \, ds
  \non
  \\
  && \leq C \delta \iot \norma{(\thetah - \thetaeps)(s)}_2 \, \normaV{(\phih - \phieps)(s)} \, ds
  \non
  \\
  \separa
  && \leq  \pier{\delta} \iot \normaV{(\phih - \phieps)(s)}^2 \, ds
  + \pier{c}\,\delta \iot \normaH{(\thetah - \thetaeps)(s)}^2 \, ds
  \non
  \\
  && \leq \pier{\delta} \iot \normaV{(\phih - \phieps)(s)}^2 \, ds
  + c \, h^2 \norma{\dt\thetaeps}_{\L2H}^2
  \leq \delta \iot \normaV{(\phih - \phieps)(s)}^2 \, ds
  + c \, h^2
  \non
\Esist
where the marked constant $C$ satisfies $\norma v_4\leq C\normaV v$ for every $v\in V$
and the last inequality is justified by \eqref{primastima}.
Next, we have
\Bsist
  && \intQt |\thetaeps| \, |\phih - \phieps|^2
  \leq c \iot \norma{\thetaeps(s)}_4 \, \norma{(\phih - \phieps)(s)}_2 \, \norma{(\phih - \phieps)(s)}_4 \, ds
  \non
  \\
  && \leq \delta \iot \normaV{(\phih - \phieps)(s)}^2 \, ds
  + c_\delta \iot \normaV{\thetaeps(s)}^2 \, \norma{(\phih - \phieps)(s)}_2^2 \, ds
  \non
\Esist
and we observe that the function $s\mapsto\normaV{\thetaeps(s)}^2$
is estimated in $L^1(0,T)$ by~\eqref{primastima}.
At this point, we choose $\delta$ small enough and apply the Gronwall lemma.
We obtain
\Bsist
  && \intQt |\thetah - \thetaeps|^2
  + \iO |\nabla( 1*\thetah - 1*\thetaeps )(t)|^2
  + \intSt |\thetahG - \thetaepsG|^2
  \non
  \\
  && \quad {}
  + \iO |(\phih - \phieps)(t)|^2
  + \intQt |\nabla(\phih - \phieps)|^2
  \non
  \\
  && \leq c \, h^2
  + c \iO \bigl( |\nabla(1*\thetaeps)(h)|^2 + |\phieps(h) - \phiz|^2 \bigr)
  + c \iG |(1*\thetaepsG)(h)|^2
  \qquad
  \label{quasisecondastima}
\Esist
for $h>0$ small enough and for every $t\in[0,T-h]$ and $\eps\in(0,1)$.
Now, we \revis{estimate} the last integrals of \eqref{quasisecondastima}
and show that they have order~$h^2$.
To this end, we \pier{argue} rather similarly as we did in deriving \eqref{persecondastima},
but we use $\thetaeps-\thetaz$ and $\phieps-\phiz$ as test functions.
Namely, we integrate \eqref{primaeps} with respect to time over $(0,s)$,
test the equality we obtain by $v=\thetaeps(s)-\thetaz$ and $\vG=\thetaepsG(s)-\thetaz\suG$,
and integrate over $(0,t)$ with respect to~$s$.
Besides, we multiply \eqref{secondaeps} by $\phieps-\phiz$ and integrate over~$Q_t$.
Finally, we add the \pier{resulting equalities to each other} and suitably rearrange.
We~have
\Bsist
  && \iO |\thetaeps(t) - \thetaz|^2
  + \frac 12 \iO |\nabla(1*\thetaeps - 1*\thetaz)(t)|^2
  + \tau \iG |\thetaepsG(t) - \thetaz\suG|^2
  \non
  \\
  && \quad {}
  + \frac 12 \iO |\phieps(t) - \phiz|^2
  + \sigma \intQt |\nabla(\phieps-\phiz)|^2
  + \intQt \bigl( \betaeps(\phieps) - \betaeps(\phiz) \bigr) (\phieps - \phiz)
  \non
  \\
  \separa
  && =
  - \intQt \bigl( \pi(\phieps) - \pi(\phiz) \bigr) (\phieps - \phiz)
  \non
  \\
  && \quad {}
  + \intQt \bigl\{
    \bigl( \thetaeps \lambdaeps(\phieps) - \thetaz \lambdaeps(\phiz) \bigr) (\phieps - \phiz)
    - \bigl( \Lambdaeps(\phieps) - \Lambdaeps(\phiz) \bigr) (\thetaeps - \thetaz)
  \bigr\}
  \non
  \\
  && \quad {}
  + \intSt m (1*u) (\thetaepsG - \thetaz\suG)
  + \intQt \nabla(1*\thetaz) \cdot \nabla(\thetaeps - \thetaz)
  \non
  \\
  && \quad {}
  + \intQt \bigl\{
    \sigma \nabla\phiz \cdot \nabla(\phieps - \phiz)
    + \bigl( - \betaeps(\phiz) - \pi(\phiz) - \thetaz \lambdaeps(\phiz) \bigr) (\phieps - \phiz)
  \bigr\} .
  \label{secondazero}
\Esist
Even though we are interested in taking $t=h$,
it is more convenient to let $t$ vary
in order to apply some Gronwall-type lemma.
So, we assume $t\in[0,h]$.
Also in this case, all the terms on the \lhs\ are nonnegative
and the first integral on the \rhs\ can be easily dealt with.
Moreover, the next term can be treated as in the above argument.
Furthermore, we have
\Bsist
  && \intSt m (1*u) (\thetaepsG - \thetaz\suG)
  \leq \intSt |\thetaepsG - \thetaz\suG|^2
  + c \intSt |1*u|^2
  \non
  \\
  && \leq  \intSt |\thetaepsG - \thetaz\suG|^2
  + c \, h^2 \norma u_2^2
  = \intSt |\thetaepsG - \thetaz\suG|^2
  + c \, h^2 \,.
  \non
\Esist
Next, we \revis{point out that}
\begin{align}
  & \intQt \nabla(1*\thetaz) \cdot \nabla(\thetaeps - \thetaz)
  \non
  \\
  & = \iO \nabla(1*\thetaz)(t) \cdot \nabla(1*\thetaeps - 1*\thetaz)(t)
  - \intQt \nabla\thetaz \cdot \nabla(1*\thetaeps - 1*\thetaz)
  \non
  \\
  & \leq \frac 14 \iO |\nabla(1*\thetaeps - 1*\thetaz)(t)|^2
  + h^2 \iO |\nabla\thetaz|^2
  + \iot \normaH{\nabla\thetaz} \, \pier{\normaH{\nabla(1*\thetaeps - 1*\thetaz)(t)}} \, dt \,.
  \non
\end{align}
Finally, the last integral of \eqref{secondazero} can be written as
\Beq
  \intQt f_\eps (\phieps - \phiz)
  \quad \hbox{where} \quad
  f_\eps := - \sigma \Delta\phiz - \betaeps(\phiz) - \pi(\phiz) - \thetaz \lambdaeps(\phiz)
  \non
\Eeq
whence it is bounded by
\Beq
  \iot \normaH{f_\eps} \, \normaH{(\phieps - \phiz)(t)} \, dt \,.
  \non
\Eeq
On the other hand, we \pier{observe that}
\Beq
  \int_0^h \normaH{\nabla\thetaz} \, dt
  = h \, \normaH{\nabla\thetaz}
  = c \, h
  \aand
  \int_0^h \normaH{f_\eps} \, dt
  = h \, \normaH{f_\eps}
  \leq c \,h
  \non
\Eeq
by \pier{virtue of} \eqref{hpdati} and~\eqref{propYosida}.
Hence, we can apply well-known Gronwall-type inequalities.
Namely we combine, e.g., \cite[Lemma~A.4, p.~156]{Brezis} and \cite[Lemma~A.5, p.~157]{Brezis}.
By ignoring some nonnegative terms on the \lhs, we conclude~that
\Beq
  \iO \bigl( |\nabla(1*\thetaeps - 1*\thetaz)(t)|^2 + |\phieps(t) - \phiz|^2 \bigr)
  + \iG |\thetaepsG(t) - \thetaz\suG|^2
  \leq c \, h^2
  \non
\Eeq
for $0\leq t \leq h$.
In particular, we have
\Bsist
  && \iO \bigl( |\nabla(1*\thetaeps)(h)|^2 + |\phieps(h) - \phiz|^2 \bigr)
  \non
  \\
  && \leq c \, h^2
  + c \iO |\nabla(1*\thetaz)(h)|^2
  \leq c \, h^2 \bigl( 1 + \normaH{\nabla\thetaz}^2 \bigr)
  = c \, h^2
  \non
\Esist
as well as
\Bsist
  && \iG |(1*\thetaepsG)(h)|^2
  = \bigl\| \textstyle \int_0^h \thetaepsG(s) \, ds \bigr\|_{\HG}^2
  \non
  \\
  && \leq 2 \bigl( \textstyle \int_0^h \norma{\thetaepsG(s) - \thetaz\suG}_{\HG} \, ds \bigr)^2
  + 2 h^2 \, \norma{\thetaz\suG}_{\HG}^2
  \leq c \, h^2 .
  \non
\Esist
Therefore, the \rhs\ of \eqref{quasisecondastima} has order~$h^2$
and the difference quotients associated to the terms of the \lhs\ of \eqref{persecondastima}
are bounded in the proper norms.
Hence, \eqref{secondastima} follows,
the bound in $\L\infty V$ for $\thetaeps$ being due to
$\nabla\thetaeps=\dt(1*\nabla\thetaeps)$.


\section*{\revis{Acknowledgments}}
This research activity has been performed in the framework of an
Italian-Romanian three-year project on ``Nonlinear partial differential 
equations (PDE) with applications in modeling cell growth, chemotaxis 
and phase transition'', financed by the Italian CNR and the Romanian 
Academy. \pcol{The present paper 
also benefits from the support of the MIUR-PRIN Grant 2010A2TFX2 ``Calculus of Variations'' and the GNAMPA (Gruppo Nazionale per l'Analisi Matematica, 
la Probabilit\`a e le loro Applicazioni) of INdAM (Istituto Nazionale di Alta 
Matematica) for PC and~GG, and of the CNCS-UEFISCDI grant, project 
PN-II-ID-PCE-2011-3-0045, for GM.}


\vspace{3truemm}

\Begin{thebibliography}{99}

\bibitem{Barbu}
V. Barbu,
\gabri{``Nonlinear differential equations of monotone types in Banach spaces'',
Springer, New York, 2010.}

\pier{
\bibitem{BBCG}
V. Barbu, M.L. Bernardi, P. Colli, G. Gilardi,
Optimal control problems of phase relaxation models,
J. Optim. Theory Appl. {\bf 109} (2001), 557-585.}

\pier{
\bibitem{BFM}
V. Barbu, A. Favini, G. Marinoschi,
Nonlinear parabolic flows with dynamic flux on the boundary, 
\pcol{J. Differential Equations {\bf 258} (2015), 2160-2195.}}

\pier{\bibitem{BCF} 
J.L. Boldrini, B.M.C. Caretta, E. Fern{\'a}ndez-Cara,
Some optimal control problems for a two-phase field model of solidification,
Rev. Mat. Complut. {\bf 23} (2010), 49-75.}

\bibitem{Brezis}
H. Brezis,
``Op\'erateurs maximaux monotones et semi-groupes de contractions
dans les espaces de Hilbert'',
North-Holland Math. Stud.
{\bf 5},
North-Holland,
Amsterdam,
1973.

\pier{\bibitem{BrokSpr} 
\pier{M. Brokate, J. Sprekels,}
``Hysteresis and Phase Transitions'',
Springer, New York, 1996.
}

\pier{
\bibitem{Cag}
G.~Caginalp,
An analysis of a phase field model of a free boundary,
Arch. Rational Mech. Anal \pier{{\bf 92} (1986), 205-245}.} 

\pier{
\bibitem{CFS}
P. Colli, M.H. Farshbaf-Shaker, J. Sprekels,
A deep quench approach to the optimal control of an AllenÐCahn equation
with dynamic boundary conditions and double obstacles,
Appl. Math. Optim. {\bf 71} (2015), 1-24.}

\bibitem{CF1}
P. Colli, T. Fukao,
The Allen-Cahn equation with dynamic boundary conditions and mass constraints,
\revis{Math.\ Methods Appl.\ Sci., DOI 10.1002/mma.3329 
(see also preprint arXiv:1405.0116 [math.AP] (2014), pp.~1-23).}

\pcol{
\bibitem{CF2} P. Colli, T. Fukao, 
Cahn--Hilliard equation with dynamic boundary conditions 
and mass constraint on the boundary, J. Math. Anal. Appl.
{\bf 429} (2015), 1190-1213.}

\bibitem{CoGiMaRo}
P. Colli, G. Gilardi, G. Marinoschi, E. Rocca,
\pier{Optimal control} for a phase field system 
with a possibly singular potential,
\pier{preprint arXiv:1410.6718~[math.AP] (2014), pp.~1-20}.

\pier{
\bibitem{CGPS}
\pier{P. Colli, G.~Gilardi,} P. Podio-Guidugli, J. Sprekels, 
Distributed optimal control of a nonstandard system of phase field equations,
Contin. Mech. Thermodyn. {\bf 24} (2012), \pier{437-459}.}
		
\pier{\bibitem{CGS} P. Colli, G. Gilardi, J. Sprekels, 
Analysis and optimal boundary control of a 
nonstandard system of phase field equations, 
Milan J. Math. {\bf 80} (2012), 119-149.}

\pcol{
\bibitem{CGS2} P. Colli, G. Gilardi, J. Sprekels, 
On the Cahn--Hilliard equation with dynamic boundary conditions
and a dominating boundary potential, J. Math. Anal. Appl. 
{\bf 419} (2014), 972-994.}

\pcol{
\bibitem{CGS3}P. Colli, G. Gilardi, J. Sprekels, 
A boundary control problem for the pure 
Cahn--Hilliard equation with dynamic boundary conditions, 
Adv.\ Nonlinear Anal., DOI 10.1515/anona-2015-0035 
(see also preprint arXiv:1503.03213 [math.AP] (2015), pp.~1-18).}

\bibitem{CMR}
P. Colli, G. Marinoschi, E. Rocca, 
Sharp interface control in a Penrose-Fife model,
\revis{ESAIM Control Optim. Calc. Var., DOI 10.1051/cocv/2015014
(see also preprint arXiv:1403.4446~[math.AP] (2014), pp.~1-33).}

\bibitem{CS}
P. Colli, J. Sprekels,
Optimal control of an Allen-Cahn equation with singular 
potentials and dynamic boundary condition,
\revis{SIAM J. Control Optim. {\bf 53} (2015), 213-234.}

\gabri{\bibitem{Ruiz}
G.R. Goldstein, 
Derivation and physical interpretation of general boundary conditions,
Adv. Differential Equations {\bf 11} (2006), 457-480.}

\gabri{\bibitem{GrMS}
M. Grasselli, A. Miranville, G. Schimperna,
The {C}aginalp phase-field system with coupled dynamic 
boundary conditions and singular potentials,
Discrete Contin. Dyn. Syst. {\bf 28} (2010), 67-98.}

\pier{\bibitem{HoffJiang}
K.-H.Hoffmann, L.S. Jiang, 
Optimal control of a phase field model for solidification,
Numer. Funct. Anal. Optim. \pier{{\bf 13} (1992), 11-27}.}

\pier{\bibitem{HKKY}
K.-H. Hoffmann, N. Kenmochi, M. Kubo, N. Yamazaki, 
Optimal control problems for models of phase-field 
type with hysteresis of play operator,
Adv. Math. Sci. Appl. {\bf 17} (2007), \pier{305-336}.}

\bibitem{LSU}
O.A. Lady\v zenskaja, V.A. Solonnikov, N.N. Ural'ceva:
``Linear and quasilinear equations of parabolic type'',
Trans. Amer. Math. Soc., {\bf 23},
Amer. Math. Soc., Providence, RI,
1968.
	
\pier{\bibitem{LeSp} C. Lefter, J. Sprekels, 
Optimal boundary control of a phase field 
system modeling nonisothermal phase transitions,
Adv. Math. Sci. Appl. {\bf 17} (2007), 181-194.}

\bibitem{Lions}
J.-L. Lions,
``\'Equations diff\'erentielles op\'erationnelles 
et probl\`emes aux limites'',
Grundlehren, Band~111,
Springer-Verlag, Berlin, 1961.

\bibitem{Lions2}
J.-L.~Lions,
``Quelques m\'ethodes de r\'esolution des probl\`emes aux limites non lin\'eaires'',
Dunod; Gauthier-Villars, Paris, 1969.

\pier{\bibitem{SY} 
K. Shirakawa, N. Yamazaki, Optimal control problems of phase field 
system with total variation functional as the interfacial energy,
Adv. Differential Equations {\bf 18} (2013), 309-350.}

\bibitem{Simon}
J. Simon,
Compact sets in the space $L^p(0,T; B)$,
Ann. Mat. Pura Appl.~(4)\/
{\bf 146} (1987), 65-96.

\pier{%
\bibitem{SprZheng}
J. Sprekels, S. Zheng, 
Optimal control problems for a thermodynamically consistent 
model of phase-field type for phase transitions, 
Adv. Math. Sci. Appl. \pier{{\bf 1} (1992), 113-125}.%
}

\gabri{\bibitem{ZLL}
J. Zheng, J. Liu, H. Liu, 
State-constrained optimal control of phase-field equations with obstacle,
Boundary Value Problems 2013, {\bf 2013}:234.}

\End{thebibliography}

\End{document}

\bye